%% file: ex_article.tex
\documentclass[]{siamart171218}

\input{ex_shared}

\begin{document}

\maketitle

\begin{abstract}
We discuss tipping phenomena (critical transitions) in nonautonomous systems using an example of a bistable ecosystem model with environmental changes represented by time-varying parameters [Scheffer et al., {\em Ecosystems}, 11 (2008), pp. 275--279]. We give simple testable criteria for the occurrence of nonautonomous tipping from the herbivore-dominating equilibrium to the plant-only equilibrium using global properties of the autonomous frozen system with fixed-in-time parameters. To begin with, we use classical autonomous bifurcation analysis to identify a codimension-three degenerate Bogdanov--Takens bifurcation: the source of a dangerous subcritical Hopf bifurcation and the organizing center for bifurcation-induced tipping (B-tipping). Then, we introduce the concept of {\em basin instability} for equilibria to identify parameter paths along which genuine nonautonomous rate-induced tipping (R-tipping) occurs without crossing any classical autonomous bifurcations. We explain nonautonomous R-tipping in terms of maximal canard trajectories and  produce nonautonomous tipping diagrams in the plane of the magnitude and rate of a parameter shift to uncover intriguing R-tipping tongues and wiggling tipping-tracking bifurcation curves. Discussion of nontrivial dynamics arising from the interaction between B-tipping and R-tipping identifies  ``points of no return” where tipping cannot be prevented by the parameter trend reversal and ``points of return tipping” where tipping is inadvertently induced by the parameter trend reversal. Our results give new insight into the sensitivity of ecosystems to the magnitudes and rates of environmental change. Finally, a comparison 
between ``tilted'' saddle-node and subcritical Hopf normal forms reveals some universal tipping properties due to basin instability, a generic dangerous bifurcation, or the combination of both.
\end{abstract}



\section{Introduction}
Tipping points are strongly nonlinear phenomena which can be described in layman's terms as 
large, sudden and often unexpected changes in the state of a system, caused by small and slow changes 
in the external inputs~\cite{scheffer2009critical,ashwin2012tipping}. The notion of a tipping point was 
popularized by Gladwell~\cite{gladwell2000tipping} and has since been used in a wide 
range of applications including climate science~\cite{lenton2008tipping,held2004detection,wieczorek2011excitability,bathiany2016beyond} 
and ecology~\cite{scheffer2009critical,scheffer2008pulse,laurance201110,boettiger2013,siteur2016ecosystems,vanselow2018,morris2002responses}. 
Scientists have identified interesting questions in relation to different tipping 
mechanisms~\cite{ashwin2012tipping,shi2016towards}, generic early warning signals near a tipping point~\cite{scheffer2009early,dakos2008slowing,scheffer2012anticipating,ritchie2016early}, 
and the possibility of preventing tipping~\cite{biggs2009turning,hughes2013living,bolt2018climate,ritchie2017inverse,alkhayuon2019} that need to be addressed in more rigorous terms. 
For example, 
Article 2 of the 1992 United Nations Framework Convention on Climate Change (UNFCCC),
which was later extended to become the
 Kyoto Protocol~\cite{protocol1997united} and the current Paris Agreement~\cite{paris2015united},  
pointed out two {\em critical factors}---the {\em level} and the {\em time frame} for changing
greenhouse gas concentrations~\cite{unfccc1992united}---suggesting that there are at 
least two tipping mechanisms of great importance to the contemporary climate. 
More generally, tipping phenomena can be classified by a type of instability and 
analyzed in more depth, although this often requires modern mathematical techniques 
beyond classical autonomous stability theory, such as geometric singular perturbations~\cite{wieczorek2011excitability,vanselow2018},
local pullback attractors~\cite{ashwin2017parameter,alkhayuon2018rate} or snapshot attractors~\cite{drotos2015,kaszas2019}, and notions of
finite-time stability~\cite{rasmussen2010,hoyer2017,duc2016}
or transient dynamics~\cite{kaszas2018}. In particular,  \cite{ashwin2017parameter} shows that much can be understood about nonautonomous tipping in one dimension from certain properties of the autonomous {\em frozen system} with fixed-in-time inputs.

In this work, we extend the discussion from~\cite{ashwin2017parameter} to a higher dimensional example from ecology~\cite{scheffer2008pulse}. 
This example exhibits a counterintuitive behavior that cannot be explained in terms of a classical autonomous bifurcation of the frozen system. The herbivore population thrives with a slow increase in the food growth rate but goes extinct
when the food growth rate increases too fast. 
From an ecological perspective, the population collapse can be discussed in terms of a ``vicious cycle"  that is inherent to various terrestrial and aquatic ecosystems, and arises from the population growth being a nonmonotone function of the food biomass~\cite{scheffer2008pulse}. 
Here, we view the population collapse as a genuine nonautonomous bifurcation.
In the spirit of  \cite{ashwin2017parameter}, we propose a framework  that uses the global property of basin boundaries in the autonomous frozen system with fixed-in-time inputs to give criteria for the occurrence of such genuine nonautonomous bifurcations in the system with time-varying inputs. This framework should be  easily accessible to applied scientists and gives new insight into nontrivial tipping phenomena  in general.

Earlier mathematical models described tipping points as classical dangerous 
bifurcations of the frozen system that occur at {\em critical levels} of an input parameter~\cite{thompson2011predicting,kuehn2011mathematical}. 
Such bifurcations have a discontinuity in the branch of exponentially stable states (attractors) at the 
bifurcation point, which explains why a system can remain near one stable state up to a 
critical level but is destined to transition to a different state past the critical 
level~\cite{thompson1994safe}. 
However, tipping points are not just classical autonomous bifurcations. 
Some systems have {\em critical rates} of parameter change, meaning that they are very 
sensitive to how fast external conditions or inputs change. Such systems can tip to a 
different state, despite the absence of any classical autonomous bifurcations, when the input parameter 
varies slowly but fast enough~\cite{wieczorek2011excitability,siteur2014beyond,scheffer2008pulse,luke2011soil,alkhayuon2019}. 
Ashwin et al. used the framework of nonautonomous dynamical systems to identify three different tipping mechanisms~\cite{ashwin2012tipping}. Bifurcation-induced tipping (B-tipping) occurs when the 
changing parameter passes through a {\em critical level} or a classical {\em dangerous bifurcation} of the frozen system, 
at which point the stable state loses stability or simply disappears. In other words,
B-tipping describes the adiabatic effects of a parameter change. Rate-induced 
tipping (R-tipping) occurs when the parameter changes faster than some {\em critical rate} 
and the system deviates from the moving stable state (attractor) sufficiently far to cross some tipping 
threshold, e.g., the boundary of the basin of attraction. In other words, R-tipping 
describes the nontrivial nonadiabatic effects of a parameter change.
Noise-induced tipping (N-tipping) occurs when noisy fluctuations drive the system past 
some tipping threshold.
Shi, Li, and Chen gave an alternative but similar classification of tipping mechanisms 
based on relative timescales of the input and of the noisy system alone~\cite{shi2016towards}. 
Additionally, tipping points can be described as either reversible or irreversible, depending 
on whether the system returns to the original stable state in the long term~\cite{wieczorek2018}.
So far, B-tipping and R-tipping have been discussed in isolation in the literature. 
However, real-world tipping phenomena  often involve a combination of different critical factors 
and different tipping mechanisms. Motivated by this observation, we 
analyze the effects of the rate of parameter change in the  ecosystem model
near the two generic dangerous bifurcations of equilibria: saddle-node and subcritical Hopf bifurcations~\cite{thompson1994safe}. Our strategy is to
\begin{itemize}
    \item introduce the concepts of {\em parameter paths} and {\em basin instability} for  equilibria in the autonomous frozen system with fixed-in-time parameters~\cite{ashwin2017parameter,wieczorek2018}  to give new insight into 
    testable criteria for genuine nonautonomous R-tipping in the system with time-varying parameters;
    \item complement classical autonomous bifurcation diagrams for the frozen system with  new information about genuine nonautonomous R-tipping bifurcations, which are entirely due to the rate of change of the input parameters,
    can be very relevant in applications, but are missed by classical autonomous bifurcation analysis;
    \item reveal nontrivial phenomena such as {\em multiple critical rates} and {\em points of no return} that arise from B-tipping due to a  dangerous autonomous bifurcation, R-tipping due to the rate of parameter change, and the interaction between the two tipping mechanisms.
\end{itemize}

Ecological models appear to be a perfect test bed for this type of study.
B-tipping  has been observed and studied extensively in different ecosystems~\cite{lewontin1969meaning,noy1975stability,scheffer1993alternative,scheffer2001catastrophic,leemans2004another,lenton2008tipping}, although the concept of a ``global tipping point" in the context of planetary 
boundaries has recently received some criticism~\cite{montoya2018planetary}. Ecologists speak of a 
``regime shift" when the bifurcation is safe or explosive and of a ``critical transition" when the 
bifurcation is dangerous~\cite{scheffer2009critical}; we refer to~\cite{thompson1994safe} for the 
classification of bifurcations into safe, explosive, and dangerous.
Similarly, there is great and rapidly growing interest in R-tipping in the context of ecological dynamics~\cite{leemans2004another,jezkova2016}.
To the best of our knowledge, the first examples of
R-tipping were reported in ecosystems~\cite{morris2002responses,scheffer2008pulse,wieczorek2011excitability,siteur2014beyond,siteur2016ecosystems,vanselow2018,neijnens2019rate}. More precisely, R-tipping conceptualizes 
a failure to adapt to changing environments~\cite{perryman2014,botero2015evolutionary}, in the sense that 
the stable state is continuously available, but the system is unable to adjust to its changing 
position when the change happens too fast. This raises the interesting research question 
of whether tipping phenomena observed 
in nature are predominantly rate-induced. What is more, the related question of whether tipping  
can be avoided or prevented has  recently received much attention in the ecosystem literature~\cite{hughes2013living,biggs2009turning,bolt2018climate,ritchie2017inverse}. Proper mathematical 
analysis of the interaction between critical levels and critical rates, or between 
B-tipping and R-tipping, is exactly what is needed to gain more insight into these questions. 
Last, there is a strong need to better understand whether ecosystems are 
sensitive to the magnitudes of environmental change, the rates of environmental change, or both. 
This is of particular importance in view of a highly variable contemporary climate, 
intensifying human activity, and rapidly declining resources.

The paper is organized as follows. Section~\ref{sec:keynon} introduces the 
ecological model given by  two nonautonomous ordinary differential equations 
and discusses the key nonlinearity. Section~\ref{sec:me} introduces the concepts 
of parameter paths and moving equilibria. 
In section~\ref{sec:Btip} we perform classical autonomous bifurcation analysis of 
the frozen system with fixed-in-time parameters and give simple criteria for
the occurrence of B-tipping in the nonautonomous system.
In section~\ref{sec:Rtip} we explain the vicious cycle that arises from the population growth being a nonmonotone function of the food biomass. We then
introduce the concept of basin instability for equilibria in the 
 autonomous frozen system to give testable criteria for  R-tipping to occur
in the nonautonomous system. We superimpose regions of basin instability on classical 
bifurcation diagrams to complement them with genuine nonautonomous R-tipping bifurcations. 
In section~\ref{sec:IRtip} we obtain two-dimensional nonautonomous R-tipping  diagrams in the parameter plane 
of the \emph{rate} and \emph{magnitude} of parameter shift and uncover R-tipping tongues.
We also describe nontrivial tipping phenomena arising from the 
interaction between B-tipping and R-tipping such as tipping diagrams with S-shaped
nonautonomous tipping-tracking bifurcation curves and multiple critical rates.
In section~\ref{sec:pnr} we partition the tipping diagrams into ``points of tracking,''
``points of return," ``points of no return," and ``points of return tipping" to give 
new insight into the problem of preventing tipping by a parameter trend reversal. 
Finally, we discuss the interaction between
B-tipping and R-tipping for modified (tilted) normal forms of the two generic dangerous bifurcations of equilibria namely saddle-node and subcritical Hopf. 
We show that the nonautonomous tipping diagram from section~\ref{sec:IRtip} 
appears to be typical for nonmonotone parameter shifts that cross
a basin instability boundary and a generic dangerous bifurcation but then reverse. 
Section~\ref{sec:concl} summarises our findings.

\section{The ecosystem model and its key nonlinearity}
\label{sec:keynon}

\begin{table}[t]
  \begin{center}
  \caption{Description of the system parameters and their values~{\normalfont\cite{scheffer2008pulse}}.}
      \begin{tabular}{cllc}
      \hline
      Symbol & Description & Units & Default value\\
      \hline
      $C > 0$ & Competition factor of plants  & m$^2$g$^{-1}$d$^{-1}$ & 0.02 \\
      $a > 0$ & Half-saturation constant of functional  & g\,m$^{-2}$ & 10 \\
      & response  &  &  \\
      $b \geq 0 $ & Exponent determining the reduced &  m$^2$g$^{-1}$ & 0--0.04\\ 
      & quality of food if food biomass  &  &  \\
      & is too high  &  &  \\
      $b_c \geq 0$ & Exponent determining the predation  & m$^2$g$^{-1}$ & 0--0.04\\ 
      &  efficiency of herbivores at high food  &  &  \\
      &  biomass  &  &  \\
      $E > 0$ & Assimilation efficiency of herbivores & dimensionless & 0.4 \\
      $c_{max} > 0$ & Maximum food intake of herbivores  & d$^{-1}$ & 1 \\ 
      & when $b_c = 0$  &  &  \\
      $m > 0$ &  Herbivore mortality rate  & d$^{-1}$ & 0--0.2 \\
      $r > 0$ & Maximum plant growth  rate& d$^{-1}$ & 0--2.5 \\
      \hline
    \end{tabular}
    \label{tab:1}
  \end{center}
\end{table}

\begin{figure}[t]
\begin{center}
\includegraphics[width=30.5pc]{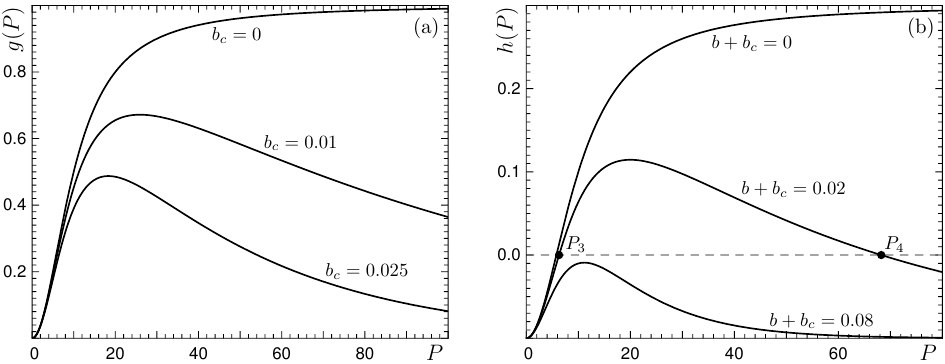}
\end{center}
\caption{{\normalfont(a)} The functional response $g(P)$ with dependence on $b_c$.
{\normalfont(b)} The key system nonlinearity: For $b + b_c > 0$, the net per-capita 
herbivore growth $h(P) = (dH/dt)/H$ has optimal plant biomass $P_{opt}$ 
where the growth is maximal and may change sign twice at $P_3$ 
and $P_4$; $m = 0.1$.}
\label{fig:graz}
\end{figure}

We consider a simple ecosystem model, where the time evolution of 
plant $P\ge 0$ and herbivore $H\ge 0$ biomass concentrations 
is modeled using two coupled nonautonomous ordinary differential 
equations~\cite{scheffer2008pulse}:
\begin{align}
\label{eq:dPdt_na}
\frac{dP}{dt} &= r(t) P - C P^2 - H\, g(P), \\
\label{eq:dHdt_na}
\frac{dH}{dt} &= \left(E\, e^{-bP} g(P)-m(t)\right)\!H,
\end{align}
with eight parameters. We fix six of the system parameters to the values or ranges 
given in Table~\ref{tab:1}. To describe changing environmental conditions, we allow 
the plant growth rate $r(t)$ and the herbivore mortality $m(t)$, which are the 
two {\em input parameters} for this study, to vary smoothly and possibly nonmonotonically over time from one asymptotic 
value to another. For example, $r(t)$ could describe the occurrence of a wet season, owing to a weather anomaly or El Niño Southern Oscillations, or changes in resources and habitat quality. Similarly, $m(t)$ could describe a disease outbreak among herbivores. 
The functional response
\begin{equation}
\label{eq:g(P)} 
g(P) = c_{max}\,\dfrac{P^2}{P^2+a^2}\,e^{-b_c P} 
\end{equation}
is a modification  of the classical monotone and strictly increasing
type-III functional response $c_{max} P^2/(P^2 + a^2)$~\cite{holling1959components} 
with an exponential factor $e^{-b_c P}$ to account for 
a decline in foraging at high plant biomass.
The resulting nonmonotone $g(P)$,
shown in Figure~\ref{fig:graz}(a) for different values of predation 
efficiency $b_c$, is believed to describe a wide range of terrestrial 
and aquatic ecosystems; see~\cite{van1996patterns,scheffer2008pulse} and 
references therein. For example, rabbits graze more with faster-growing 
plants as long as the plants are small enough but avoid overgrown 
bushes for fear of predators and are unable to graze on 
plants that have grown too tall. Similarly, in aquatic ecosystems, 
phytoplankton can be heavily consumed at early life stages by 
herbivorous zooplankton, but higher-density phytoplankton colonies 
become less prone to exploration and foraging.
%
%
Our aim is to give criteria for tipping in the nonautonomous 
system~\eqref{eq:dPdt_na}--\eqref{eq:dHdt_na} with time-varying $r(t)$ or $m(t)$ 
in terms of certain local and global properties of the autonomous {\em frozen system}
\begin{align}
\label{eq:dP/dt}
\frac{dP}{dt} &= r P - C P^2 -  H\, g(P), \\
\label{eq:dH/dt}
\frac{dH}{dt} &= \left(E\, e^{-bP} g(P)-m\right)\!H, 
\end{align}
where $r$ and $m$ are fixed-in-time input parameters. 

Owing to the modified functional response~\eqref{eq:g(P)}, 
the frozen system~\eqref{eq:dP/dt}--\eqref{eq:dH/dt} 
is a singular perturbation problem: it has a different number 
of equilibrium solutions for $b+b_c =0$ and  $0<b+b_c\ll 1$. 
To see that, consider the net per-capita herbivore growth $h(P) =(dH/dt)/H$, 
shown in Figure~\ref{fig:graz}(b),
whose roots correspond to nonzero herbivore equilibrium concentrations.
When $b + b_c = 0$, the net per-capita herbivore growth
is a strictly increasing function of $P$ with a single root $P_3$ 
(Figure~\ref{fig:graz}(b)). However, when $0<b+b_c\ll 1$, the net per-capita 
herbivore growth has a maximum at the optimal plant biomass
$$
P_{opt} \approx  \left(\frac{2a^2}{b + b_c}\right)^{\!\frac{1}{3}}
$$
and can have no roots at all, one double root, or two distinct roots 
at $P_3 < P_{opt}$  and $P_4 > P_{opt}$ (Figure~\ref{fig:graz}(b)); see~\cite{OKeeffePhD2020} 
for the derivation of $P_{opt}$, $P_3$, and $P_4$.
This key nonlinearity underpins nonautonomous R-tipping bifurcations and
arises from the decline in foraging at high 
plant biomass $(b_c > 0)$, from reduced food quality at high plant 
biomass $(b > 0)$, or from a combination of both (Figure~\ref{fig:graz}(b)).
Thus, we refer to $b+ b_c$ as the 
{\em nonlinearity parameter} and work with  different but fixed-in-time 
values of $b$ and $b_c$, as indicated in Table~\ref{tab:1}.

\section{Moving equilibria and parameter paths}
\label{sec:me}

Typically, the position of
an equilibrium $e$ for the frozen system depends on the input parameters $r$ and/or $m$.
When the input parameters vary over time, $e$ changes its position
in the $(P,H)$ phase space
and we speak of a {\em moving equilibrium} 
$$
e(t) := e(r(t),m(t)),
$$
also known as a quasistatic equilibrium~\cite{ashwin2012tipping}.
Note that $e(t)$ is a property of the autonomous 
frozen system~\eqref{eq:dP/dt}--\eqref{eq:dH/dt} and the changing 
environment, but it is not a solution to the nonautonomous 
system~\eqref{eq:dPdt_na}--\eqref{eq:dHdt_na}.

As the input parameters $r(t)$ and $m(t)$ evolve smoothly over time,  
they trace out a continuous {\em parameter path} $\Delta$ in the two-dimensional 
$(r,m)$ parameter plane:
$$
\Delta := \left\{(r(t),m(t)): t\in\mathbb{R}   \right\}.
$$
We use the notions of a  {moving stable equilibrium}
and a {parameter path} to discuss the differences and interactions 
between B-tipping and R-tipping.

\section{B-tipping: Classical autonomous bifurcations}
\label{sec:Btip}

Before we analyze genuine nonautonomous R-tipping bifurcations in the original nonautonomous system \eqref{eq:dPdt_na}--\eqref{eq:dHdt_na}, we perform classical 
autonomous bifurcation analysis of the frozen system~\eqref{eq:dP/dt}--\eqref{eq:dH/dt}~\cite{kuznetsovelements}.
This will allow us to contrast two tipping mechanisms: B-tipping, which is due to a slow parameter drift via a classical dangerous bifurcation of the frozen system, and  genuine nonautonomous R-tipping, which is entirely due to the rate of parameter change and need not involve any classical autonomous bifurcations of the frozen system.
Specifically, we compute bifurcation curves in the $(r,m)$ parameter plane and uncover 
a  codimension-three degenerate Bogdanov--Takens bifurcation: the source of a dangerous subcritical Hopf bifurcation and the organizing center for B-tipping. We say an 
equilibrium or a limit cycle is stable when it is exponentially stable.

\subsection{Steady bifurcations of equilibria}

The frozen system~\eqref{eq:dP/dt}--\eqref{eq:dH/dt} has at most four equilibria. 
A {\em trivial} equilibrium $e_1$ and a \emph{plant-only} equilibrium $e_2$,
$$
e_1  = (0,0),\quad e_2  = \left(r/C,0\right),
$$
exist for all parameter settings. The trivial equilibrium $e_1$ is a saddle with eigenvalues
$\lambda_1 = r > 0$ and $\lambda_2 = - m<0$.
The plant-only equilibrium $e_2$ has eigenvalues $\lambda_1 = -r<0$ and 
$\lambda_2 = E c_{max}\, e^{-(b+b_c)r/C}/\left((a\,C/r)^2 + 1\right) -m$.
Hence, $e_2$ is a stable node when $\lambda_2 < 0$, is a saddle when 
$\lambda_2 > 0$, and undergoes a transcritical or pitchfork bifurcation
whenever $\lambda_2 = 0$.

The other two equilibria correspond to a stationary coexistence of plants and 
herbivores and satisfy the following conditions:
\begin{align}
\label{eq:nonzeroHequilibrium2}
H &= \frac{(r - CP)(P^2 + a^2)}{c_{max} \,P \,e^{-b_c P}},\\
\label{eq:nonzeroHequilibrium1}
h(P) &= E\,c_{max} \,\frac{P^2 e^{-(b+b_c) P}}{P^2+a^2} - m = 0.
\end{align}
Although the roots of  \eqref{eq:nonzeroHequilibrium1} cannot be expressed in a closed form, one can take advantage of the small nonlinearity parameter 
$0< b+b_c\ll 1$ and use perturbation methods to obtain closed form 
approximations in terms of an asymptotic expansion in different powers of 
$b+b_c$; see~\cite[sect.3.1]{OKeeffePhD2020} for the details of the derivations.
Regular perturbation about $b + b_c = 0$ gives the $P$-\,component of 
the \emph{herbivore-dominating} equilibrium $e_3$:
\begin{equation}
\label{eq:pert1} 
P_3 = \sqrt{\frac{a^2\, m}{E\,c_{max}-m}} + \frac{a^2\, m\,E\,c_{max}}{2\left(E\,c_{max} - m\right)^{\,2}}\, (b+b_c) + \mathcal{O}\!\left((b+b_c)^2\right), 
\end{equation}
where $\mathcal{O}((b+b_c)^2)$ is the error term of order 
$(b+b_c)^2$ as $(b+b_c)\to 0$, and
\begin{equation} 
\label{eq:e3} 
e_3 = \left(\sqrt{\dfrac{a^2\, m}{E\,c_{max}-m}} +  \mathcal{O}(b+b_c),\frac{(r - CP_3)(P_3^2 + a^2)}{c_{max}\, P_3 \,e^{-b_c P_3}}\right).
\end{equation}
Singular perturbation about $b + b_c = 0$ using a stretched variable 
$\tilde{P} = (b + b_c) P$ gives the $P$-\,component of 
the \emph{plant-dominating} equilibrium $e_4$:
\begin{equation}
\label{eq:pert2} 
P_4 = \frac{\ln\!\left(E\,c_{max}/m\right)}{b+b_c} -  
\frac{a^2(b+b_c)}{\left(\ln\!\left(E\,c_{max}/m\right)\!\right)^{\,2}}+ \mathcal{O}\!\left((b+b_c)^2\right), 
\end{equation}
and
\begin{equation}
\label{eq:e4} 
e_4 = \left(\frac{\ln\!\left(E\,c_{max}/m\right)}{b+b_c} +  \mathcal{O}(b+b_c),\frac{(r - CP_4)(P_4^2 + a^2)}{c_{max}\,P_4 \,e^{-b_c P_4}}\right). 
\end{equation}
The main advantage of the closed form approximations is the 
information about the dependence of the equilibrium positions on the 
system parameters $r$, $m$ and $(b+b_c)$. 

Next, we consider the qubic equation
\begin{equation}
\label{eq:7} 
q(P) = (b+b_c)P^3 + a^2(b+b_c)P - 2\,a^2 = 0,
\end{equation}
for the parameter values from Table~\ref{tab:1}, and we show as follows.
\begin{proposition}
In the $(r,m)$ parameter plane of the frozen system~\eqref{eq:dP/dt}--\eqref{eq:dH/dt}, there is a transcritical bifurcation curve 
\begin{equation}
\label{eq:trans1} 
T = \left\{ (r,m): r\in\mathbb{R}_+\setminus \{CP^*\},\,
m = \frac{E\,c_{max} e^{-(b+b_c)r/C}}{(a\,C/r)^2 + 1}
\right\}.
\end{equation}
If  \eqref{eq:7} has a root $P^*>0$, then there is a half-line of saddle-node bifurcations
\begin{equation}
\label{eq:saddle1} 
S_e = \left\{(r,m): r > C P^*,\;m = \dfrac{E\,c_{max}\,e^{-(b + b_c) P^*}}{(a/P^*)^2 + 1} \right\}
\end{equation}
and a pitchfork bifurcation point
$$
Pf = \left\{(r,m): r = C P^*,\; m = \frac{E\,c_{max} e^{-(b+b_c)P^*}}{(a/P^*)^2 + 1}\right\}.
$$
\end{proposition}

{\em Proof of Proposition 4.1.}
Equilibrium  $e_3$ or $e_4$  becomes degenerate with equilibrium $e_2$ in a 
transverse crossing if $P=r/C$ in 
Eqs.~\eqref{eq:nonzeroHequilibrium2}--\eqref{eq:nonzeroHequilibrium1}. 
The crossing corresponds to a codimension-one transcritical bifurcation or to a codimension-two (due to the lack of the $\mathbb{Z}_2$-symmetry) pitchfork bifurcation.
Thus, substituting $P=r/C$ into  \eqref{eq:nonzeroHequilibrium1} 
defines the curve $T$ of transcritical bifurcations in the $(r,m)$ parameter plane that may include pitchfork bifurcation points.
Equilibria $e_3$ and $e_4$ become degenerate in a quadratic (saddle-node) or cubic
(pitchfork) tangency when $r$-independent  \eqref{eq:nonzeroHequilibrium1} has a 
positive repeated root, meaning that 
\begin{equation}
\label{eq:12} 
h(P) = 0\; \text{ and }\; h'(P) = 0\; 
\text{ for some}\; P > 0.
\end{equation}

One can verify that the second equation in~\eqref{eq:12} holds if and only if the cubic equation~\eqref{eq:7} holds. Since $q'(P) \ge 0$ and $q(0) < 0, h'(P)$ can have at most one positive root $P=P^*$. This
root is used in~\eqref{eq:nonzeroHequilibrium1} to determine the value of $m$ at which $e_3$ and $e_4$ become degenerate.
To eliminate a triple degeneracy involving $e_1$ or $e_2$, we require that the corresponding $H$
from~\eqref{eq:nonzeroHequilibrium2} is positive, meaning that $r > CP^*$.
Thus, conditions~\eqref{eq:12} together with $r>CP$ define the half-line $S_e$
of saddle-node bifurcations of equilibria in the $(r,m)$ parameter 
plane.
Finally, note that curves $T$ and $S_e$ meet at the point $Pf$, which corresponds to a 
cubic tangency of $h(P)$ or a
triple degeneracy of $e_2$, $e_3$, and $e_4$. This is a pitchfork bifurcation point.
$\blacksquare $

If $b + b_c = 0$ (Figure~\ref{fig:rm}(a)), then  \eqref{eq:7} has no roots, 
meaning that there are no saddle-node or pitchfork bifurcations 
of equilibria. The curve $T$ of transcritical bifurcations emerges from the origin 
and levels off at $m = E\,c_{max}$ for large $r$. 
The equilibrium $e_3$ that bifurcates from $e_2$ along $T$ exists
below $T$.  
If $b + b_c > 0$ (Figure~\ref{fig:rm}(b)--(d)), then  \eqref{eq:7} has a unique root 
$P^*>0$ that corresponds to a unique half-line $S_e$ of saddle-node bifurcations. 
Equilibria $e_3$ and $e_4$ that are born along $S_e$ exist below $S_e$.
The curve $T$ of transcritical bifurcations emerges from the origin, 
has a maximum $Pf$, and approaches $m=0$ from above 
for large $r$. Now, $T$ consists 
of two different branches separated by $Pf$.
Equilibrium $e_3$, which bifurcates from $e_2$ along the 
solid branch of $T$, exists below the solid branch of $T$.
In contrast, equilibrium $e_4$, which bifurcates from $e_2$ along the 
dashed branch of $T$, exists above the dashed branch of $T$.
Equilibria $e_2$, $e_3$, and $e_4$ become degenerate in a pitchfork 
bifurcation at $Pf$.

\begin{figure}[t]
\begin{center}
\includegraphics[width=30.5pc]{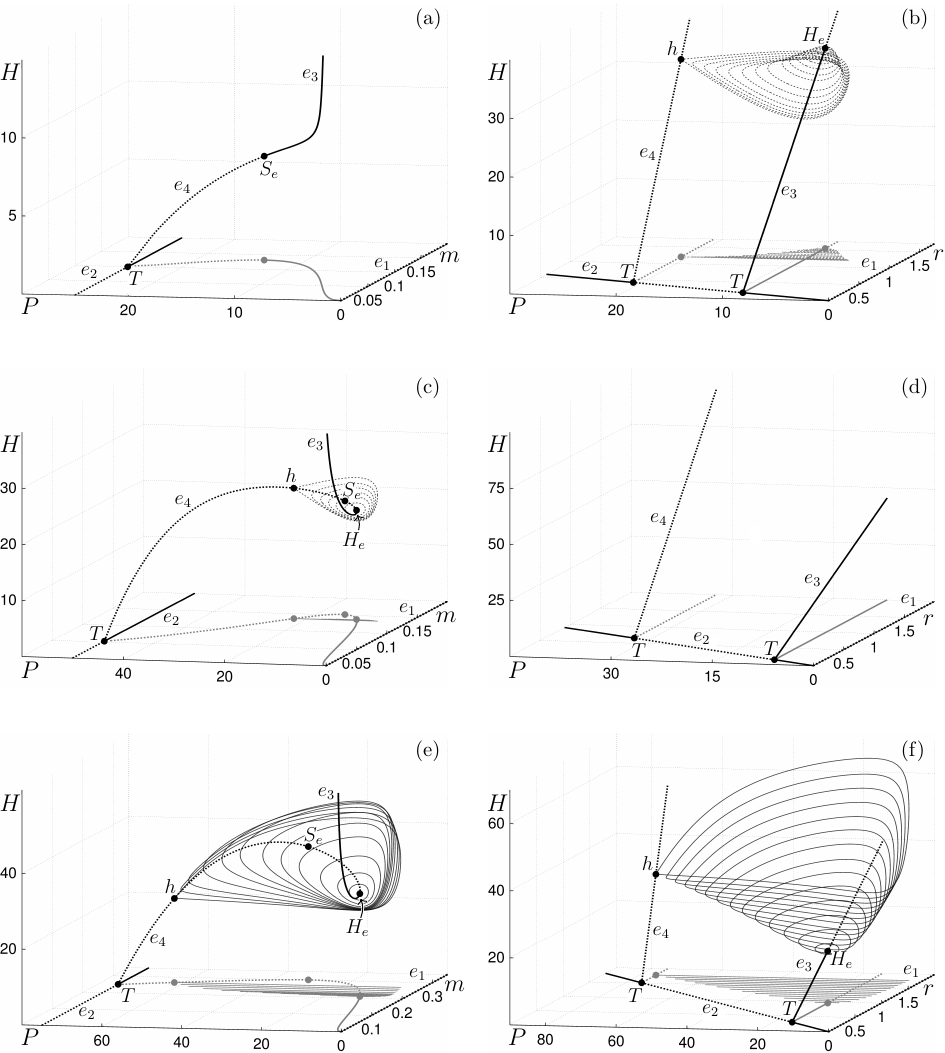}
\end{center}
\caption{One-parameter autonomous bifurcation diagrams for the 
frozen system~\eqref{eq:dP/dt}--\eqref{eq:dH/dt}
showing the position and stability 
of equilibria and limit cycles. 
The left column shows  the $(P,H,m)$ space for {\normalfont(a)} $r=0.5$, $(b,b_c)= (0.025,0.025)$, 
{\normalfont(c)} $r=1$, $(b,b_c)= (0.02,0.02)$, and {\normalfont(e)} $r=1.5$, $(b,b_c)= (0.001,0.005)$.
The right column shows the $(P,H,r)$  space for {\normalfont(b)} $m=0.115$, 
$(b,b_c)= (0.025,0.025)$, {\normalfont(d)} $m=0.1$, $(b,b_c)= (0.02,0.02)$, and {\normalfont(f)} $m=0.25$, 
$(b,b_c)= (0.001,0.005)$. Solid  branches indicate stable 
solutions, dashed branches indicate unstable solutions.  Projections onto 
the $(m,P)$ and $(r,P)$ planes are shown in gray. For the labeling of different 
bifurcations see Table~{\normalfont\ref{tab:3}}.
}
\label{PH}
\end{figure}

\subsection{Bifurcations of limit cycles}

To reveal bifurcations of limit cycles and showcase different types of autonomous 
dynamics in the autonomous frozen system~\eqref{eq:dP/dt}--\eqref{eq:dH/dt}, 
we plot six examples of one-dimensional bifurcation diagrams in Figure~\ref{PH} 
for two types of parameter paths. In the left column we fix $r$ and 
consider a range of $m\in(0,0.2]$.
In the right column we fix $m$ and consider a range of $r\in(0,2]$. 
In addition to the transcritical $T$ and saddle-node $S_e$ bifurcations
of equilibria identified in the previous section, there are 
supercritical (Figure~\ref{PH}(e)--(f)) and dangerous subcritical (Figure~\ref{PH}(b)--(c))
Hopf bifurcations $H_e$. Additionally, a limit cycle can connect to the saddle equilibrium 
$e_4$ and disappear in a homoclinic bifurcation $h$ (Figure~\ref{PH}(b)--(c) 
and (e)--(f)). Finally, there are saddle-node bifurcations of limit cycles discussed
in the next section. For more details and background on classical autonomous 
bifurcation theory, we refer to~\cite{kuznetsovelements}.

\begin{table}[t]
  \begin{center}
  \caption{Glossary of terms for bifurcation diagrams.}
    \begin{tabular}{cl}
      \hline
      Symbol & Description\\
      \hline
       $T$ & Transcritical bifurcation \\ 
       $S_e$ & Saddle-node of equilibria bifurcation \\
       $Pf$ & Pitchfork bifurcation \\
      $H_e$ & Hopf bifurcation \\ 
      $h$ & Homoclinic bifurcation \\ 
     $BT_{I(II)}$ & Bogdanov--Takens type-I(II) bifurcation \\
      $GH$ & Generalized Hopf (Bautin) bifurcation \\ 
     $S_{lc}$ & Saddle-node of limit cycles bifurcation \\
      $h_{res}$ & Resonant homoclinic bifurcation \\ 
     $BI$ & Basin instability \\ 
     \hline
    \end{tabular}
    \label{tab:3}
  \end{center}
\end{table}

\begin{figure}[t]
\begin{center}
\includegraphics[width=30.5pc]{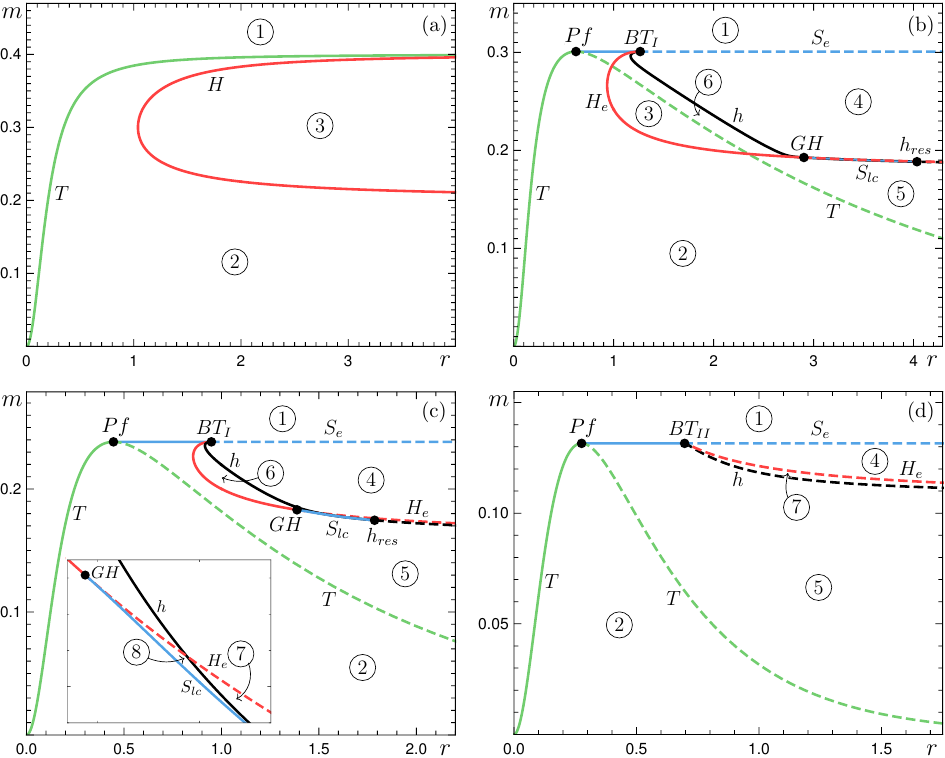}
\caption{Two-parameter autonomous bifurcation diagrams
for the frozen system~\eqref{eq:dP/dt}--\eqref{eq:dH/dt} in the $(r,m)$
parameter plane, obtained 
for different but fixed $(b, b_c)$ = {\normalfont(a)} $(0,0)$, {\normalfont(b)} $(0.001,0.005)$, {\normalfont(c)} $(0.005,0.01)$, {\normalfont(d)} $(0.025,0.025)$. 
Bifurcations that give rise to attractors are plotted as solid 
curves. For the labeling of different bifurcations see Table~{\normalfont\ref{tab:3}}.
}
\label{fig:rm}
\end{center}
\end{figure}

\subsection{Two-parameter autonomous bifurcation diagrams}

To provide a systematic bifurcation analysis, we obtain 
two-dimensional $(r,m)$ bifurcation diagrams for different 
but fixed values of $b$ and $b_c$ (Figure~\ref{fig:rm}). 
We plot codimension-one  bifurcations that give rise to attractors as
solid curves. To be more specific:
along a solid (dashed) transcritical bifurcation, 
the bifurcating branch of equilibria is stable 
(of saddle type); along a solid (dashed) saddle-node 
bifurcation, a saddle collides with an attractor 
(repeller); and along  solid (dashed) Hopf 
and homoclinic bifurcations, the bifurcating limit 
cycle is attracting (repelling). Transcritical and 
saddle-node bifurcations of equilibria are obtained 
using conditions~\eqref{eq:trans1} and~\eqref{eq:saddle1}, 
respectively. Hopf, homoclinic, and saddle-node bifurcations 
of limit cycles are computed using the numerical continuation 
software AUTO~\cite{doedel2007auto}.

\begin{figure}[t]
\begin{center}
\includegraphics[width=30.5pc]{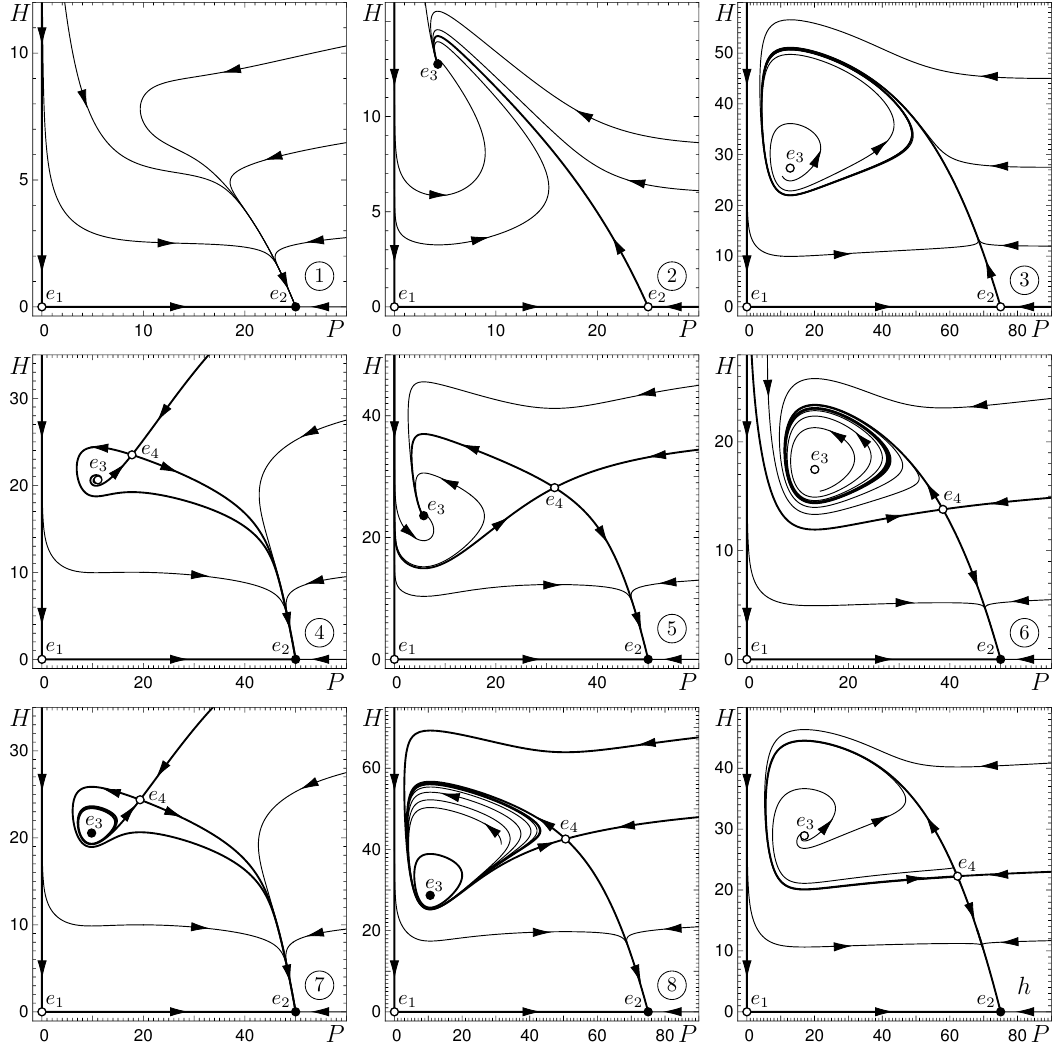}
\end{center}
\caption{Examples of qualitatively different $(P,H)$  phase portraits 
for the autonomous frozen system~\eqref{eq:dP/dt}--\eqref{eq:dH/dt} showing 
(filled circles) stable equilibria, (open circles) unstable equilibria, 
(thick curves) limit cycles and stable/unstable invariant 
manifolds of saddle equilibria, and (thin curves) examples of 
typical trajectories. Note the stable limit cycle in regions 
3 and 6, the unstable limit cycle in region 7, and the two 
limit cycles in region 8. See Table~{\normalfont\ref{tab:2}} for 
parameter values.}
\label{fig:pp}
\end{figure}

For each bifurcation diagram, we identify regions with
qualitatively different dynamics and illustrate these 
with examples of phase portraits in the $(P,H)$ phase plane 
(Figure~\ref{fig:pp}). 
It turns out that there are at least four qualitatively different 
$(r,m)$ bifurcation diagrams, depending on the 
settings of $b$ and $b_c$.

In the absence of the key nonlinearity, that is, when $b + b_c = 0$
and $g(P)$ is the classical monotone type-III functional response, 
there are just two bifurcation curves:
(solid) curve $T$ of  transcritical bifurcations and (solid) curve $H_e$ of
supercritical Hopf bifurcations (Figure~\ref{fig:rm}(a)). 
These two curves do not interact, and they partition the $(r,m)$ parameter 
plane into three distinct regions with qualitatively different 
dynamics (Figure~\ref{fig:pp}, {1}--{3}). In particular, 
$H_e$ gives rise to a stable limit cycle in region {3}, which 
represents a stable but oscillatory coexistence between plants 
and herbivores. These simple dynamics change drastically in the 
presence of the key nonlinearity.

When $b + b_c$ becomes small but nonzero, meaning the functional response $g(P)$ becomes nonmonotone, a number of 
qualitative changes take place in the autonomous bifurcation diagram as
expected from the singular perturbation nature of the problem.
Specifically, there are three additional codimension-one bifurcation curves, and four 
special codimension-two bifurcation points (Figure~\ref{fig:rm}(b)). 
Understanding the new bifurcation diagram is reminiscent of 
assembling a jigsaw-puzzle. 
First, a half-line $S_e$ of saddle-node bifurcations 
of equilibria appears. $S_e$ emerges from the 
pitchfork bifurcation point $Pf$ on $T$, 
where $T$ changes from solid to dashed. 
Second, $H_e$ is no longer unbounded at both ends
but emerges from the 
Bogdanov--Takens bifurcation point $BT_I$ on $S_e$, 
where $S_e$ changes from solid to dashed. There are 
two possible types of Bogdanov--Takens bifurcation, 
and $BT_I$ is type-I according to 
the classification in~\cite[sect.8.4]{kuznetsovelements}. 
It is known from the unfolding of a Bogdanov--Takens 
bifurcation that the curve of homoclinic bifurcations $h$
must emerge from $BT_I$. Along $h$, the limit cycle 
originating from $H_e$ becomes a connecting orbit to the saddle 
equilibrium $e_4$ and disappears (Figure~\ref{fig:pp}, {$h$}).
Third, there is a generalized Hopf (or Bautin) bifurcation point 
$GH$ on $H_e$, where $H_e$ changes from (solid) supercritical to 
(dashed) subcritical~\cite[sect.8.3]{kuznetsovelements}. 
It is known from the unfolding of a generalized Hopf 
bifurcation that 
the curve of a saddle-node of a limit cycles $S_{lc}$
must emerge from $GH$. Along solid $S_{lc}$, two limit cycles 
of which one is attracting and the other repelling collide and disappear. 
Finally, $S_{lc}$ terminates on $h$  at a resonant homoclinic
bifurcation point $h_{res}$, where $h$ changes from solid 
to dashed. This new bifurcation structure has five additional 
regions {4}--{8} with qualitatively different dynamics.

When the combination of $b$ and $b_c$ is increased further, 
 points $GH$ and $h_{res}$ approach 
$BT_I$ (Figure~\ref{fig:rm}(c)). 
In the process, region {3} with stable self-sustained 
oscillations disappears, while the bistable region {5} 
becomes noticeably larger. 
Then, there are {special combinations} of $b$ and $b_c$, 
where $GH$ and $h_{res}$ collide simultaneously with $BT_I$ 
and disappear in a codimension-three degenerate Bogdanov--Takens 
bifurcation (not shown in the figure)~\cite[sect.3.2.1]{OKeeffePhD2020}. This collision 
eliminates $S_{lc}$ together with the (dashed) supercritical part of $H_e$ 
and with regions {6} and {8}. Concurrently, 
the Bogdanov--Takens bifurcation point changes to type-II.
The difference from $BT_{I}$ is that $H_e$ and $h$ emerging 
from $BT_{II}$ swap their relative positions and become (dashed) subcritical
(Figure~\ref{fig:rm}(d)). 

\begin{table}[t]
  \begin{center}
  \caption{Parameter values chosen for phase portraits in Figure~{\normalfont\ref{fig:pp}}}
    \begin{tabular}{cllll}
      \hline
      Phase portrait & $r$ & $m$ & $b$ & $b_c$\\
      \hline
      {1} & 0.5  & 0.14 & 0.025 & 0.025 \\ 
      {2} & 0.5  & 0.05 & 0.025 & 0.025 \\ 
      {3} & 1.5  & 0.23 & 0.001 & 0.005 \\ 
      {4} & 1    & 0.125 & 0.025 & 0.025 \\ 
      {5} & 1    & 0.075 & 0.025 & 0.025 \\ 
      {6} & 1    & 0.21 & 0.005 & 0.01 \\ 
      {7} & 1    & 0.12 & 0.025 & 0.025 \\ 
      {8} & 1.5  & 0.18025 & 0.005 & 0.01 \\ 
      {\small{$h$}} & 1.5  & 0.2684 & 0.001 & 0.005 \\ 
      \hline
    \end{tabular}
    \label{tab:2}
  \end{center}
\end{table}

Past the special combination of $b$ and $b_c$, 
there are four bifurcation curves, 
including the two dangerous bifurcations of equilibria that 
are of interest for B-tipping: 
the (solid) half-line $S_e$ of saddle-node bifurcations, 
and the (dashed) curve $H_e$ of subcritical Hopf bifurcations.
Now, there are two special bifurcation points: the pitchfork point $Pf$ and a type-II Bogdanov--Takens bifurcation point $BT_{II}$.
$H_e$ gives rise to a repelling limit cycle in region {7},
which becomes a connecting orbit to the saddle equilibrium $e_4$ and disappears in a 
homoclinic bifurcation along $h$. 
Finally, a substantial part of 
the diagram is occupied by adjacent regions {5} 
and {7} with bistability between the plant-only equilibrium $e_2$ and 
the herbivore-dominating equilibrium $e_3$. This bistability
is of interest for R-tipping from $e_3$ to $e_2$ studied in section~\ref{sec:Rtip}.

\subsection{Testable criterion for B-tipping}

If a continuous parameter 
path $\Delta$ in the $(r,m)$ bifurcation  diagram crosses 
a
dangerous bifurcation for the autonomous frozen system~\eqref{eq:dP/dt}--\eqref{eq:dH/dt}, then there is a 
time-varying external input 
$(r(t),m(t))$ 
that traces out this path 
and gives rise to B-tipping in the nonautonomous system~\eqref{eq:dPdt_na}--\eqref{eq:dHdt_na}.

\begin{figure}[]
\begin{center}
\hspace*{-0.5cm}
\includegraphics[width=31pc]{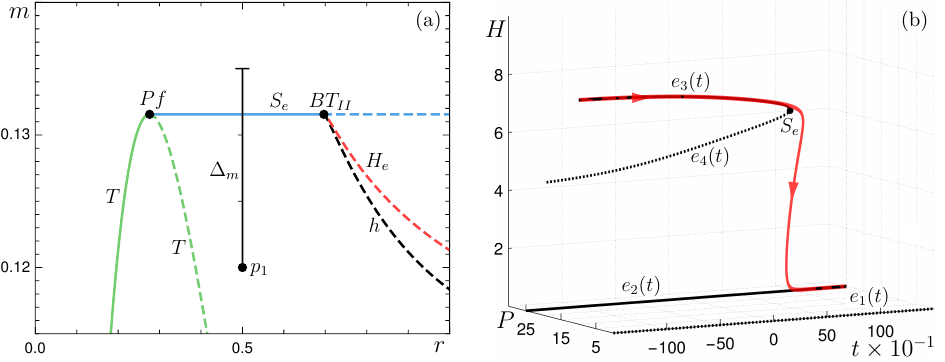}
\end{center}
\caption{{\normalfont(a)} Example of a parameter path $\Delta_m$, in the autonomous 
$(r,m)$  bifurcation diagram of the frozen system~\eqref{eq:dP/dt}--\eqref{eq:dH/dt},
that crosses a (dangerous) saddle-node 
bifurcation $S_e$. 
{\normalfont(b)} As $m(t)$ is increased from $p_1 = (0.5,0.12)$ along the path, 
the nonautonomous system~\eqref{eq:dPdt_na}--\eqref{eq:dHdt_na} 
undergoes B-tipping from $e_3(t)$ to $e_2(t)$ as $m(t)$ passes through $S_e$.
$b = b_c = 0.025$, and $m(t) = 0.12 + 0.015 (\tanh(\varepsilon t) + 1)/2$
with $\varepsilon = 10^{-3}$.
}
\label{fig:pht0}
\end{figure}

When $b + b_c = 0$, we do not expect any B-tipping owing to the lack 
of  dangerous bifurcations.  However, when $b + b_c > 0$, meaning that 
there is a decline in herbivore growth at high plant biomass, a number 
of different B-tipping mechanisms appear in the ecosystem model. 
The most dominant are the two generic dangerous bifurcations 
of equilibria, namely saddle-node and 
subcritical Hopf bifurcations. 
Figure~\ref{fig:pht0} shows an example of a parameter path, denoted $\Delta_m$ in panel (a), that crosses a dangerous saddle-node bifurcation 
$S_e$ of the frozen system~\eqref{eq:dP/dt}--\eqref{eq:dH/dt},
together with the dynamics of the 
nonautonomous system~\eqref{eq:dPdt_na}--\eqref{eq:dHdt_na}
where $m(t)$ drifts slowly along the path (panel (b)). 
If the system starts near the stable 
equilibrium $e_3$ at the lower end $p_1$ of the path and 
 $m(t)$ increases over time, then the nonautonomous 
system tracks the moving stable equilibrium $e_3(t)$ up to the point 
of the dangerous bifurcation $S_e$, which defines the 
{\it critical level} of $m$ (Figure~\ref{fig:pht0}(b)). 
As $m(t)$ passes through the bifurcation, the system 
undergoes a sudden and 
abrupt transition to the other stable equilibrium $e_2(t)$.
This transition, called here {\em B-tipping}, is also known as 
a {\em dynamic} or {\em adiabatic 
bifurcation}~\cite{Benoit1991,Rasmussen2007}.

\begin{figure}[]
\begin{center}
\hspace*{-0.5cm}
\includegraphics[width=31pc]{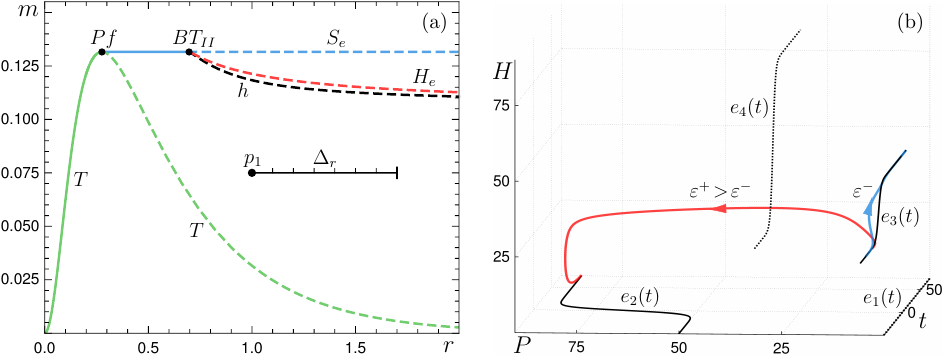}
\end{center}
\caption{{\normalfont(a)} Example of a parameter path $\Delta_r$, 
in the autonomous 
$(r,m)$  bifurcation diagram of the frozen system~\eqref{eq:dP/dt}--\eqref{eq:dH/dt},
that does not cross 
any autonomous bifurcations. 
{\normalfont(b)} As $r(t)$ 
is increased from $p_1 = (0.75,0.075)$ along the path at a rate 
$\varepsilon^-$ (blue trajectory) the nonautonomous system~\eqref{eq:dPdt_na}--\eqref{eq:dHdt_na} tracks 
the moving stable equilibrium $e_3(t)$. However, for a faster 
rate $\varepsilon^+ > \varepsilon^-$ (red trajectory) there is
irreversible R-tipping from $e_3(t)$ to $e_2(t)$ even though $e_3(t)$ 
never disappears or loses stability in the autonomous sense.
$b = b_c = 0.025$, and $r(t) = 0.75 + 0.6 (\tanh(\varepsilon t) + 1)/2$, 
with $\varepsilon^- = 0.1$ and $\varepsilon^+ = 0.2$. The moving equilibria 
are obtained for $\varepsilon \approx 0.14$.
}
\label{fig:pht01}
\end{figure} 

\section{Nonautonomous R-tipping: Beyond classical autonomous bifurcations}
\label{sec:Rtip}

In this section we go beyond the classical autonomous bifurcation theory and 
adiabatic effects of a parameter change. Specifically, we consider genuine 
nonautonomous bifurcations and nontrivial nonadiabatic effects of a 
parameter change
that arise solely from the rate of change of the input parameters $r$ and $m$
and cannot be captured by classical autonomous bifurcation analysis.
Specifically, we ask,
{\em Are there parameter paths in the autonomous $(r,m)$ bifurcation diagram 
that do not cross any  bifurcation of the stable equilibrium $e_3$ but give rise to tipping from $e_3$ to $e_2$
when the input parameter $r$ or $m$ varies over time?}
The 
answer is yes. This was demonstrated in~\cite{scheffer2008pulse} 
and is examined here in more depth. Consider a parameter path $\Delta_r$ in Figure~\ref{fig:pht01}(a) that does not cross any classical autonomous bifurcations. If the nonautonomous system starts 
at the stable equilibrium $e_3$ near the lower end  $p_1$  of the path,
and $r(t)$ increases slowly enough along the path, then the 
nonautonomous system is able to adapt to the changing environment
and track the moving stable equilibrium $e_3(t)$ 
(blue trajectory in Figure~\ref{fig:pht01}(b)). 
However, if $r(t) $ increases slowly but faster than some critical rate, 
the nonautonomous system fails to adapt to the changing environment and 
undergoes a critical transition from $e_3(t)$ to the other stable 
equilibrium $e_2(t)$ (red trajectory in Figure~\ref{fig:pht01}(b)).
Tipping occurs even though $e_3(t)$ is continuously available and never 
loses stability along the path in the autonomous sense.
Such a genuine nonautonomous bifurcation is known as irreversible {\em R-tipping}~\cite{ashwin2012tipping,wieczorek2018}.

\subsection{The vicious cycle}
\label{sec:vc}

Intuitively, R-tipping is the result of
a {\em vicious cycle} that could potentially tip 
the system to a different state if the input parameters vary 
too fast. In the ecosystem model, the vicious cycle arises 
from the key nonlinearity, namely nonmonotone 
herbivore growth $h(P) =(dH/dt)/H$ that changes sign from positive to
negative at high plant biomass $P=P_4$ (see Figure~\ref{fig:graz}(b)). 

The effect can be understood as follows.
Consider a stable herbivore population with a lower than optimal 
plant biomass $P_3$ for some $r=r_-$. Then, consider a smooth increase 
in the plant growth rate from $r_-$ to $r_+$. This results 
in faster-growing plants and moves the stable equilibrium to a 
larger herbivore population with the same plant biomass $P_3$.
If $r(t)$ increases slowly enough, herbivores manage to graze 
and grow fast enough so that the larger herbivore population 
is able to maintain the same plant biomass at larger $r=r_+$. However, 
if $r(t)$ increases too fast, herbivores may be unable to keep 
up and prevent the plant biomass from growing past its optimal 
value $P_{opt}$. This, in turn, triggers the vicious cycle: past the optimal 
plant biomass, the larger the plant biomass gets, the less the 
herbivores graze and grow, allowing the plant biomass to grow even larger. 
The ultimate effect is negative net 
herbivore growth causing a sudden collapse of the herbivore 
population. This is accompanied by a sudden increase in the plant 
biomass to $P_4$. There is no classical autonomous bifurcation along the
parameter path between $r_-$ and $r_+$, but the rate of change of $r(t)$ alone 
prevents the system from adapting to the modified stable equilibrium. 
We now introduce the key mathematical concepts to analyze
the vicious cycle mechanism that gives rise to genuine nonautonomous
R-tipping bifurcations.

\subsection{Basin instability on a path}
\label{sec:tcTtip}

\begin{figure}[t]
\begin{center}
\hspace*{-0.5cm}
\includegraphics[width=31pc]{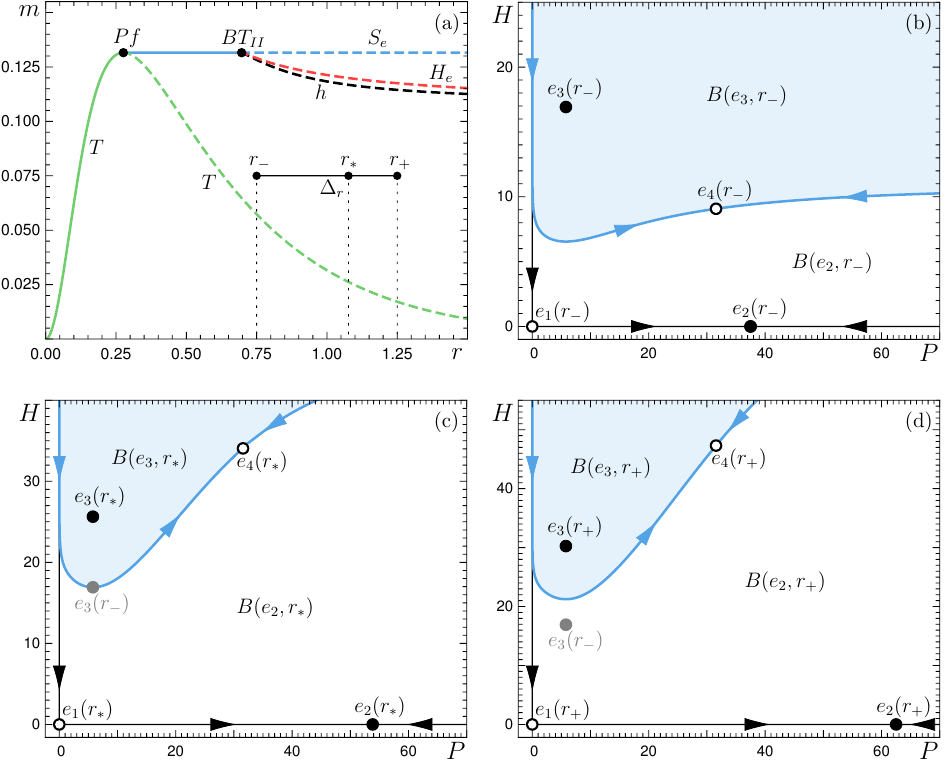}
\end{center}
\caption{{\normalfont(a)} A 
two-parameter autonomous bifurcation diagram
for the frozen system~\eqref{eq:dP/dt}--\eqref{eq:dH/dt} in the $(r,m)$
parameter plane with a parameter path $\Delta_r$ that does not cross any autonomous bifurcation. 
{\normalfont(b)--(d)} Phase portraits of~\eqref{eq:dP/dt}--\eqref{eq:dH/dt} at three different
points along the path $\Delta_r$ illustrate basin instability 
of $e_3$ on $\Delta_r$. Blue shading indicates the basin of  
attraction of stable equilibrium $e_3$ for {\normalfont(b)} $ r= r_- = 0.75$, 
{\normalfont(c)} $r = r_* \approx 1.07672$, and {\normalfont(d)} $r = r_+ = 1.25$. 
$b = b_c = 0.025$.
}
\label{fig:fbs0}
\end{figure}

It turns out that, similarly to B-tipping, much can be understood about 
genuine nonautonomous 
R-tipping in system~\eqref{eq:dPdt_na}--\eqref{eq:dHdt_na} 
from certain properties of the autonomous 
frozen system~\eqref{eq:dP/dt}--\eqref{eq:dH/dt}~\cite{ashwin2017parameter}. 
The difference is
that R-tipping is related to global, rather than local, 
properties of the stable equilibrium. 
The key concept for understanding 
irreversible R-tipping is {\em basin instability}, and 
we need the following ingredients  to define it:
\begin{itemize}
\item[(i)] 
An {\em exponentially stable base equilibrium} $e(p)$ of the frozen system, whose position 
in the phase space varies with the input parameter(s) $p$. In the ecosystem model, 
the stable equilibrium $e_3$ is given in terms of 
$p=(r,m)$ 
by  \eqref{eq:e3}.
\item[(ii)]
{\em Bistability or multistability} in the frozen system: at 
least one attractor in addition to $e(p)$ for the same input parameters. 
The ecosystem model exhibits bistability between $e_3$ and $e_2$
in the $(r,m)$ parameter regions 5 and 7.
\item[(iii)] A continuous {\em parameter path} $\Delta$ in the 
autonomous bifurcation diagram of the frozen system that does not cross any 
dangerous bifurcation of a stable base equilibrium $e(p)$. For example, see the path $\Delta_r$ in the $(r,m)$ bifurcation diagram of the ecosystem model
in Figure~\ref{fig:fbs0}(a).
%
\item[(iv)] {\em The basin of attraction} of a stable base equilibrium $e(p)$
in the frozen system, defined as the set of initial states $(P_0,H_0)$ whose trajectories 
converge to $e(p)$:
$$
B(e,p) = \{ 
(P_0,H_0)\in\mathbb{R}^2: (P(t),H(t))\to e(p),\; t\to +\infty
\}.
$$
\end{itemize}
%
%
\begin{definition}
In the autonomous frozen system~\eqref{eq:dP/dt}--\eqref{eq:dH/dt}, consider a continuous parameter path $\Delta$ with a family (branch) of exponentially stable equilibria $e(p)$ that vary $C^1$-smoothly
with $p\in\Delta$. Let $\overline{B(e,p)}$ denote the basin of attraction of $e(p)$ together with its boundary.
We say that $e(p)$ is {\em basin unstable on the path 
$\Delta$} if there are two points on the path $p_1, p_2\in\Delta$, such that $e(p_1)$ 
is outside the basin of attraction of $e(p_2)$:
$$
e(p_1) \notin \overline{B(e,p_2)}.
$$
\end{definition}
Note that basin instability is a global property of the autonomous frozen system and a chosen parameter path.
This is illustrated for the herbivore-dominating equilibrium $e_3(r)$ on the path $\Delta_r$  in Figure~\ref{fig:fbs0}. The stable equilibrium $e_3(r_-)$ is contained
within the basin of attraction of $e_3(r)$ for $r\in[r_-, r_*)$, 
lies on the basin boundary of $e_3(r_*)$ (Figure~\ref{fig:fbs0}(c)), 
and is outside the basin of attraction of $e_3(r)$ for
$r\in(r_*,r_+]$ (Figure~\ref{fig:fbs0}(d)).  Thus, $e_3(r)$ is 
basin unstable on $\Delta_r$.
To include the global basin instability property in the classical autonomous bifurcation 
diagram, we make the following definition.
\begin{definition}
In the autonomous frozen system~\eqref{eq:dP/dt}--\eqref{eq:dH/dt}, consider $e(p_1)$ and $e(p_2)$ from a $C^1$-smooth family of exponentially stable equilibria $e(p)$. We define the region of basin instability of $e(p_1)$ as the set of points $p_2$ in the parameter space such that $e(p_1)$ lies outside the basin of attraction of $e(p_2)$:
\begin{align}
\label{eq:BI}
BI(e,p_1) = \{p_2: e(p_1)\notin \overline{B(e,p_2)}\}.
\end{align}
\end{definition}
The region of basin instability of the herbivore-dominating equilibrium $e_3(p_1)$ in parts $5$ and $7$ of the $(r,m)$  plane, denoted $BI$ and defined as
\begin{align}
\label{eq:BIe}
BI := BI(e_3,p_1) = \{p_2\in\circled{5}\cup\circled{7}: e_3(p_1)\notin \overline{B(e_3,p_2)}\},
\end{align}
is shown in gray in Figure~\ref{fig:pht1}(a).

We refer to~\cite{wieczorek2018} for extension 
of these ideas to {\em threshold instability} to also capture reversible R-tipping that does not require bistability or basin boundaries.\footnote{In the transient 
phenomenon of reversible R-tipping, the system fails to track the moving stable state and suddenly 
moves to a different state, but in the long term returns to and tracks
the original stable state~\cite{wieczorek2011excitability,vanselow2018, wieczorek2018}.}

\subsection{Testable criterion for R-tipping and maximal canards} 

Consider a continuous parameter path $\Delta$.
Suppose there is a family (branch) of exponentially stable equilibria $e(p)$ of the autonomous frozen system~\eqref{eq:dP/dt}--\eqref{eq:dH/dt}
that vary $C^1$-smoothly with $p\in\Delta$, meaning there is no classical autonomous bifurcation of $e(p)$ on $\Delta$.
If  $e(p)$ is basin unstable on $\Delta$, meaning that there are points $p_1,p_2\in\Delta$ such that 
$e(p_1) \notin \overline{B(e,p_2)}$, then there 
is a time-varying external input 
$p(t) = (r(t),m(t))$ that traces out 
the path from $p_1$ to $p_2$ and gives irreversible
R-tipping from $e(t)$ in the nonautonomous system~\eqref{eq:dPdt_na}--\eqref{eq:dHdt_na}~\cite{wieczorek2018}.

It can be shown rigorously that basin instability is necessary and sufficient 
for the occurrence of R-tipping in one-dimensional systems~\cite{ashwin2017parameter} 
and sufficient but not necessary for the occurrence of R-tipping in higher-dimensional systems~\cite{wieczorek2018,kiers2018conditions,xie2020}. Here, we explain the rigorous results intuitively,
using the example of a parameter path $\Delta_r$ from Figure~\ref{fig:fbs0}.
Suppose the nonautonomous system is initialized 
in the basin of attraction and near the stable equilibrium 
$e_3(r_-)$, then undergoes a monotone parameter shift from 
$r_-$ to $r_+$.
If $r(t)$ varies sufficiently slowly, 
the nonautonomous system is guaranteed to closely track (adiabatically follow) 
the moving stable equilibrium $e_3(t)$ along the 
path~\cite{ashwin2017parameter,wieczorek2018}. 
If  $r(t)$ shifts abruptly, then just after the shift, the nonautonomous system remains
at its earlier position near $e_3(r_-)$.  This now lies outside 
the basin of attraction of $e_3(r_+)$ and inside the basin of 
attraction of $e_2(r_+)$ (Figure~\ref{fig:fbs0}(d)), so the system
converges to $e_2(r_+)$. Thus, there must be at least one intermediate critical 
rate of change of $r(t)$ at which the nonautonomous tracking-tipping bifurcation occurs.

\begin{figure}[t]
\begin{center}
\hspace*{-0.5cm}
\includegraphics[width=31pc]{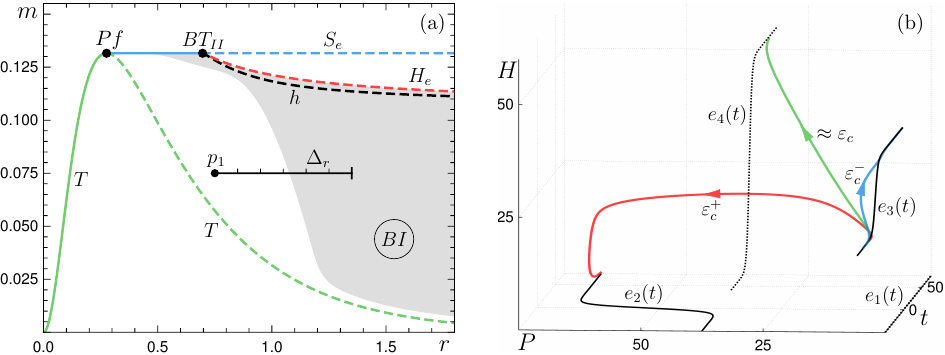}
\end{center}
\caption{
The same as Figure~{\normalfont\ref{fig:pht01}} but with {\normalfont(a)} the extended 
path $\Delta_r$ and the addition of the shaded region of basin 
instability $BI := BI(e_3,p_1)$ for 
$p_1 = (0.75,0.075)$ as defined by  \eqref{eq:BIe}, and {\normalfont(b)}
the addition of the green canard trajectory for 
$\varepsilon \approx \varepsilon_c$ that somewhat surprisingly
tracks the moving unstable equilibrium $e_4(t)$.
}
\label{fig:pht1}
\end{figure}

Analysis of the region of basin instability of $e_3(p_1)$ in Figure~\ref{fig:pht1}(a) reveals that genuine nonautonomous R-tipping bifurcations are ubiquitous in the ecosystem model. They will occur on every parameter path that connects $p_1$ to some $p_2\in BI$ and stays within regions 5 and 7. Thus,
in addition to dangerous magnitudes of environmental change, the ecosystem model 
appears to be particularly sensitive to how fast the plant growth rate $r$ increases 
over time.  The proposed  concept of basin instability  quantifies this rate-of-change 
sensitivity and can be applied to any nonlinear system. Superimposing the region of 
basin instability onto a classical autonomous bifurcation diagram gives basic information 
on genuine nonautonomous bifurcations that can be very relevant in certain applications but 
are missed by classical autonomous bifurcation analysis. 

One may ask about the dynamics at a critical rate, where a transition between tracking and R-tipping occurs. It turns out that, at a critical 
rate, the nonautonomous system somewhat surprisingly tracks the moving 
unstable equilibrium $e_4(t)$ for an infinite time (Figure~\ref{fig:pht1}(b)). 
In the terminology of slow-fast systems, a genuine nonautonomous R-tipping 
bifurcation corresponds to a unique {\em maximal canard trajectory} that remains 
within an unstable slow manifold for the longest time\footnote{
To see that, reformulate the two-dimensional 
nonautonomous system~\eqref{eq:dPdt_na}--\eqref{eq:dHdt_na}
as a three-dimensional slow-fast autonomous system  by 
augmenting it with $u=\varepsilon t$ as an additional dependent 
variable so that $du/dt =\varepsilon$. Then,  in the $(P,H,u)$ phase space, $e_4(u)$ becomes 
a normally hyperbolic unstable critical manifold for $\varepsilon=0$, 
which persists as a normally hyperbolic unstable slow manifold for 
$0 < \varepsilon \ll 1$.}~\cite{szmolyan2001canards}. 
Depending on the basin boundary in the
frozen system, critical-rate trajectories can track other moving unstable states 
such as limit cycles~\cite[Figure 5.1]{OKeeffePhD2020}, which are referred to as 
{\em edge states} in~\cite{wieczorek2018}.
Alternatively, a critical-rate trajectory can be transformed into a connecting (heteroclinic) orbit using the compactification technique developed in~\cite{wieczorek2020}, as shown in~\cite[Chap.7]{OKeeffePhD2020} and~\cite{xie2020}.

\section{Nonautonomous tipping diagrams for parameter shifts}
\label{sec:IRtip}

Guided by the R-tipping criterion and basin instability analysis performed 
in the previous section, we analyze the nonautonomous system~\eqref{eq:dPdt_na}--\eqref{eq:dHdt_na} 
with a monotone shift 
\begin{equation}
\label{eq:force1} 
r(t) = r_- + \dfrac{\Delta_{r}}{2}\left(\tanh(\varepsilon t) + 1\right),
\end{equation}
from $r_-$ to $r_+ = r_- + \Delta_r$ with $\dot{r}_{\scriptsize{max}} = \varepsilon \Delta_r/2$, 
and a nonmonotone shift 
\begin{equation}
\label{eq:force2} 
r(t) = r_- + \Delta_{r} \, \mbox{sech}(\eps t), 
\end{equation}
from $r_-$ to $r_+ = r_- + \Delta_r$ and then back to $r_-$ with
$\dot{r}_{\scriptsize{max}} = \varepsilon \Delta_r/2$. 
The shifts are parameterized by their {\em magnitude} 
$\Delta_r$ and {\em rate} $\varepsilon$,
which
enables parametric study in the form of two-dimensional $(\Delta_r, \varepsilon)$ 
or $(\Delta_r, \dot{r}_{max})$ {\em nonautonomous tipping diagrams}. In this way, we  identify 
{\em critical rates} $\varepsilon_c$ at which the system undergoes a 
nonautonomous bifurcation from tracking to (irreversible) R-tipping.

\subsection{Monotone shifts across a basin instability boundary: Single critical rate}

\begin{figure}[t]
\begin{center}
\hspace*{-0.5cm}
\includegraphics[width=31pc]{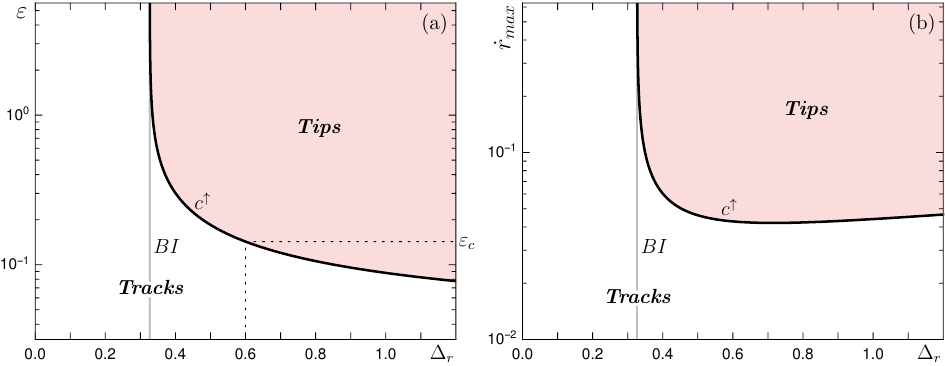}
\end{center}
\caption{Nonautonomous tipping diagrams in the {\normalfont(a)} $(\Delta_r,\varepsilon)$ 
and {\normalfont(b)} $(\Delta_r,\dot{r}_{max})$ parameter planes for monotone 
shifts~\eqref{eq:force1} from $p_1 = (0.75,0.075)$  along the 
extended parameter path $\Delta_r$ from Figure~{\normalfont\ref{fig:pht1}(a)}. 
The nonautonomous tipping-tracking bifurcation curve $c^\uparrow$ separates the 
diagram into regions of (white) tracking and (pink) irreversible 
R-tipping. The vertical gray line  indicates the boundary of the 
basin instability region $BI := BI(e_3,p_1)$ defined by  \eqref{eq:BIe}.
The critical rate $\varepsilon_c$ corresponds 
to the (green) canard trajectory in Figure~{\normalfont\ref{fig:pht1}(b)}. 
$b = b_c = 0.025$.
}
\label{fig:epsrdot1}
\end{figure} 

\begin{figure}[t]
\begin{center}
\hspace*{-0.5cm}
\includegraphics[width=31pc]{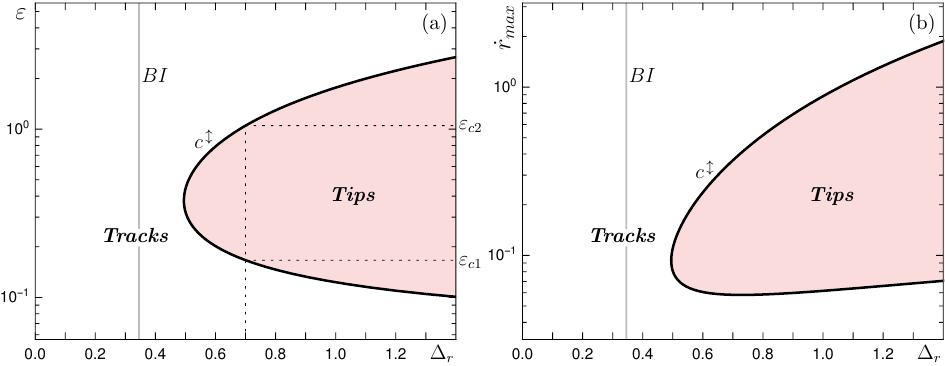}
\end{center}
\caption{
Nonautonomous tipping diagrams in the {\normalfont(a)} $(\Delta_r,\varepsilon)$ 
and {\normalfont(b)} $(\Delta_r,\dot{r}_{max})$ parameter plane for 
nonmonotone shifts~\eqref{eq:force2} from $p_1 = (1.0,0.075)$ 
along a path $\Delta_r$ with a fixed $m=0.075$ and varied $r>1$.
The nonautonomous tipping-tracking bifurcation curve $c^{\,\updownarrow}$ separates the 
diagram into regions of (white) tracking and (pink) irreversible 
R-tipping. The vertical gray line  indicates the boundary of the basin instability region 
$BI := BI(e_3,p_1)$ defined by  \eqref{eq:BIe}. $b = b_c = 0.025$.
}
\label{fig:epsrdot2}
\end{figure}  


A systematic analysis of R-tipping for monotone shifts~\eqref{eq:force1} 
from $p_1$ along the path $\Delta_r$ from Figure~\ref{fig:pht1}(a) gives the 
$(\Delta_r,\varepsilon)$ and $(\Delta_r,\dot{r}_{max})$
tipping diagrams in Figure~\ref{fig:epsrdot1}. 
The nonautonomous R-tipping bifurcations 
occur along the curve $c^\uparrow$, which 
 divides the tipping diagram into separate regions of (white) tracking and 
(pink) irreversible R-tipping (Figure~\ref{fig:epsrdot1}).  
As $\Delta_r$ decreases, the $c^\uparrow$ curve becomes asymptotic to the (gray line) 
boundary of the basin instability region $BI$.
Furthermore, the $c^\uparrow$ curve appears to level off at $\dot{r}_{max}\approx 0.045$.
Thus, one can give simple approximate conditions for the occurrence of
irreversible R-tipping along this path in terms of the shift magnitude $\Delta_r$ exceeding the boundary
of the basin instability $BI$, and $\dot{r}_{max}$ exceeding the critical value $\approx 0.045$.
Finally, the monotone shift has a unique critical rate 
$\varepsilon_c$ for a fixed magnitude $\Delta_r$.

\subsection{Nonmonotone shifts across a basin instability boundary: Two critical rates}

Analysis of R-tipping for nonmonotone shifts~\eqref{eq:force2} 
tracing out the path $\Delta_r$ in
Figure~\ref{fig:pht1}(a) from $p_1$ at $r_- = 0.75 $ to 
$r_- + \Delta_r$ and then back to $p_1$ is shown in the tipping diagram in Figure~\ref{fig:epsrdot2}.
The nonautonomous R-tipping bifurcation curve $c^{\,\updownarrow}$ 
forms an {\em R-tipping tongue} which is reminiscent of a resonance tongue for 
time-periodic inputs~\cite{marchionne2018synchronisation},
in the sense that the system exhibits a strongly enhanced response to external 
inputs with optimal timing. 
As $\varepsilon$ is decreased from above, 
the natural timescales of $H(t)$ and $P(t)$ get closer to the timescale of $e_3(t)$, the 
system starts to react to the impulse input in $r$, and R-tips due to basin instability. 
This transition is marked by the higher critical rate $\varepsilon_{c1}$. As
the natural timescales of $H(t)$ and $P(t)$ become 
comparable to the timescale of $e_3(t)$, there is a strongly enhanced 
response in the form of 
a tipping tongue. As $\varepsilon$ is decreased even further, 
the natural timescales of $H(t)$ and $P(t)$ become faster than the timescale of $e_3(t)$, 
and the system starts to closely track $e_3(t)$. This transition 
is marked by the lower critical rate $\varepsilon_{c1}$. The nonmonotone shift across a basin instability boundary
 typically has two critical rates, $\varepsilon_{c1}$ and $\varepsilon_{c2}$,  
for a fixed magnitude $\Delta_r$.

\begin{figure}[t]
\begin{center}
\hspace*{-0.5cm}
\includegraphics[width=31pc]{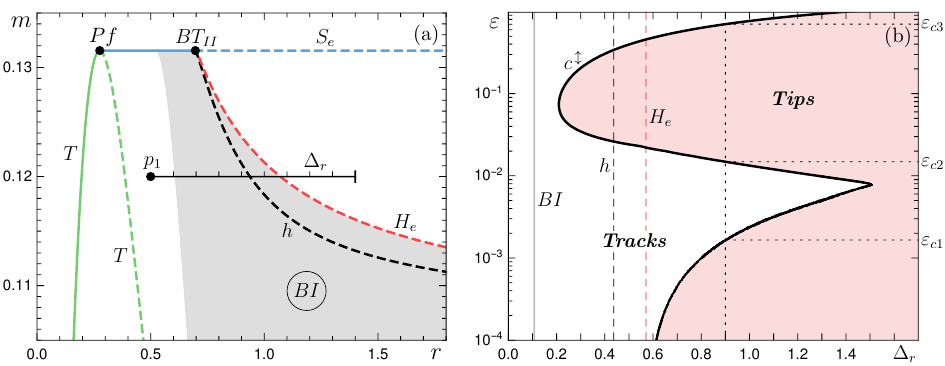}
\end{center}
\caption{
{\normalfont(a)} 
Example of a parameter path $\Delta_r$, in the autonomous 
$(r,m)$  bifurcation diagram of the frozen system~\eqref{eq:dP/dt}--\eqref{eq:dH/dt},
that crosses the boundary of the basin instability region $BI := BI(e_3,p_1)$
for $p_1 = (0.5,0.12)$ as defined by  \eqref{eq:BIe},
and the (dangerous) subcritical Hopf bifurcation $H_e$.
{\normalfont(b)} The nonautonomous tipping diagram in the $(\Delta_m,\varepsilon)$ parameter 
plane for nonmonotone shift~\eqref{eq:force2} from $p_1$ along 
the path $\Delta_r$ from panel {\normalfont(a)}. The nonautonomous tipping-tracking bifurcation 
curve $c^{\,\updownarrow}$ separates the diagram into regions of (white) 
tracking and (pink) tipping.
$b = b_c = 0.025$.
}
\label{fig:pht4}
\end{figure}

\subsection{Nonmonotone shifts across a basin instability boundary and a dangerous bifurcation: A critical level and three critical rates}

So far, we have discussed B-tipping and R-tipping in isolation. This and the following 
subsections reveal interesting tipping phenomena that arise from the interaction between B-tipping and R-tipping, or between critical levels and critical rates. 

Consider nonmonotone shifts along the path $\Delta_r$ in Figure~\ref{fig:pht4}(a) 
from $p_1$, past the boundary of the  basin instability $BI$, past the (dangerous) subcritical 
Hopf bifurcation $H_e$, and back to $p_1$.
Now, the nonautonomous tipping-tracking
bifurcation curve $c^{\,\updownarrow}$  consists of two distinct parts,
which correspond to two different tipping mechanisms (Figure~\ref{fig:pht4}(b)). 
 At high $\varepsilon$ and 
between $BI$ and $H_e$, we replicate the distinctive R-tipping tongue from Figure~\ref{fig:epsrdot2} and attribute this 
part to pure irreversible R-tipping.
As $\varepsilon$ is decreased, there are two new features.
First, the curve $c^{\,\updownarrow}$ forms a deep  wedge
whose tip delineates the change from R-tipping to B-tipping. Second, 
as $\varepsilon\to 0$, the curve $c^{\,\updownarrow}$ approaches 
the critical level $H_e$ for B-tipping, but the approach is very
``slow''. The new features can be explained in terms of relative timescales 
and a bifurcation delay. As $\varepsilon$ decreases below the tipping tongue, 
the natural timescales of $H(t)$ and $P(t)$ 
start to exceed the timescale of $e_3(t)$, meaning that the system becomes more able to 
follow the moving stable equilibrium $e_3(t)$. On the one hand, we 
start to lose R-tipping. On the other hand, the system acquires some 
characteristics of a slow passage through a Hopf bifurcation. 
As $\varepsilon$ is decreased further, $H(t)$ and $P(t)$ become much 
faster than $e_3(t)$ and start to closely track $e_3(t)$. We move into the regime 
of a slow passage through a Hopf bifurcation, which is characterised
by a noticeable bifurcation delay that does not vanish even when the 
rate of parameter change tends to zero~\cite{baer1989slow,neishtadt1987persistence,neishtadt1988persistence}.
This means that trajectories follow the moving equilibrium past the bifurcation 
point, where the equilibrium becomes unstable, for a noticeable time even when $\varepsilon\to 0$.
Thus, the deep tracking wedge in $c^{\,\updownarrow}$ and the ``slow'' approach of $c^{\,\updownarrow}$ 
toward $H_e$ as $\varepsilon\to 0$ are attributed to this bifurcation delay. 
Finally, the  change in the basin boundary at $h$,
from  the stable invariant manifold of the saddle $e_4$ to the unstable limit cycle, may 
be another contributing factor.
In summary, the intricate tipping diagram captures 
different aspects of the interaction between B-tipping and R-tipping, 
giving rise to three critical rates.

\begin{figure}[t]
\begin{center}
\hspace*{-0.5cm}
\includegraphics[width=31pc]{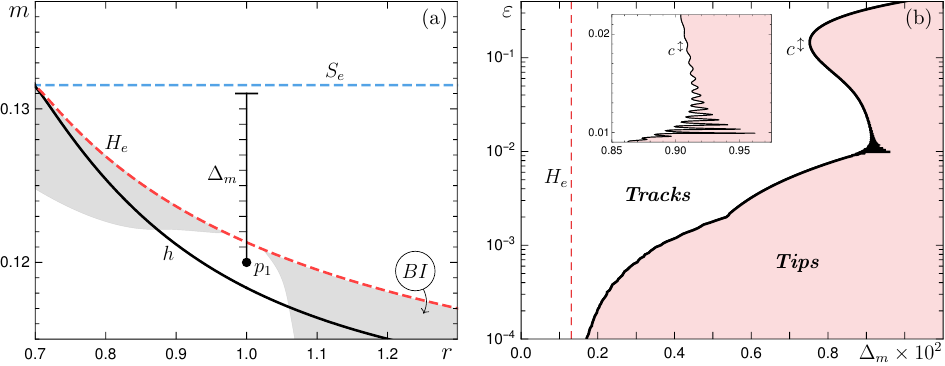}
\end{center}
\caption{
{\normalfont(a)} Example of a parameter path $\Delta_m$, in the autonomous 
$(r,m)$  bifurcation diagram of the frozen system~\eqref{eq:dP/dt}--\eqref{eq:dH/dt},
that crosses the dangerous subcritical Hopf bifurcation $H_e$.
(Gray) The basin instability region $BI := BI(e_3,p_1)$ for $p_1= (1,0.12)$ is
defined by  \eqref{eq:BIe}.
{\normalfont(b)} The nonautonomous tipping diagram in the $(\Delta_m,\varepsilon)$ parameter 
plane for nonmonotone shift~\eqref{eq:force2} from $p_1$ along 
the path $\Delta_m$ from panel {\normalfont(a)}. The nonautonomous tipping-tracking bifurcation 
curve $c^{\,\updownarrow}$ separates the diagram into regions of (white) 
tracking and (pink) tipping. 
$b = b_c = 0.025$.
}
\label{fig:fbs10}
\end{figure}

\subsection{Nonmonotone shifts across a dangerous bifurcation: A critical level and multiple critical rates}

To reveal an interesting tipping effect that arises near a (dangerous) 
subcritical Hopf bifurcation, we fix $r = 1$ and consider nonmonotone shifts in the 
herbivore death rate,
$
m(t) = 0.12 + \Delta_m \sinh(\varepsilon t),
$
along a path $\Delta_m$ from Figure~\ref{fig:fbs10}(a) that crosses a (dangerous) 
subcritical Hopf bifurcation $H_e$ with a vanishing region of basin instability.

The resulting nonautonomous tipping-tracking bifurcation curve $c^{\,\updownarrow}$ shows a complicated 
rate dependence that is far from trivial (Figure~\ref{fig:fbs10}(b)). Owing to 
the absence of basin instability, no pure R-tipping occurs along this path. 
Nonetheless, there can be  multiple critical rates. 
Past $H_e$, there is a range of shift magnitudes $\Delta_m$ 
with a unique critical rate. However, for larger shift magnitudes, 
the curve $c^{\,\updownarrow}$ has a remnant of the R-tipping tongue
that gives rise to three critical rates for a fixed $\Delta_m$. Most interestingly, 
there is an interval of $\Delta_m$ 
where the wiggling part of $c^{\,\updownarrow\,}$ gives rise to 
several critical rates for a fixed $\Delta_m$ (inset in Figure~\ref{fig:fbs10}(b)).

\section{Points of return, points of no return, points of return tipping}
\label{sec:pnr}

Tipping is often defined as a large, sudden, and possibly unexpected change 
in the state of the system, caused by a slow or small change in the external 
input (e.g., environmental conditions). Although ``sudden" and ``unexpected" 
suggest that foreseeing and preventing tipping may be difficult, it should 
in general be possible~\cite{hughes2013living}. In this section, we are 
guided by the question: 
{\em Given a monotone parameter shift that gives tipping, under what 
conditions can tipping be 
prevented by a parameter-shift reversal?} Certain aspects of this question 
have been explored in the context of B-tipping near a saddle-node bifurcation. 
For example, Hughes et al.~\cite{hughes2013living} speak of ``living dangerously on 
borrowed time'' to describe a window of opportunity for ecosystems to 
return to safer conditions before an otherwise inevitable tipping occurs.
Biggs et al.~\cite{biggs2009turning} ask whether early-warning indicators for 
tipping provide sufficient warning to modify the ecosystem's management and 
avert undesired regime shifts by ``turning back from the brink."
Gandhi, Knobloch, and Beaume~\cite{gandhi2015dynamics,gandhi2015localized} consider nonmonotone 
parameter shifts through the (global) saddle-node on a limit cycle bifurcation 
to identify a new resonance mechanism in 
the context of spatially localized (vegetation) patterns.
Ritchie, Karabacak, and Sieber~\cite{ritchie2017inverse} model systems near a saddle-node bifurcation 
and analyze the relationship between the time and amplitude of a saddle-node 
crossing to avoid B-tipping. A similar problem is analyzed by Li et al.~\cite{li2019time} in terms of pullback attractors and points of no return.
More recently, Alkhayuon et al. investigate ``avoided" B-tipping
and R-tipping near a subcritical Hopf bifurcation in the box model of 
the Atlantic Meridional Overturning Circulation (AMOC) in the context of 
the collapse of the AMOC
and climate change mitigation~\cite{alkhayuon2019}.

Here, we extend the existing literature on avoiding B-tipping to
include R-tipping effects due to basin instability.
Specifically, we consider  a nonautonomous system with paths in one parameter $\mu$
that may but do not need to cross a dangerous bifurcation at $\mu=\mu_b$. 
Along a  path, we keep the nonmonotone shift~\eqref{eq:force1} unchanged and 
modify the monotone shift~\eqref{eq:force2} to reach a maximum in finite time: 
\begin{equation}
  \mu(t) = \begin{cases}
     \mu_- + \Delta_{ \mu} \, \mbox{sech}(\eps t), &  \,  t \leq 0, \\
     \mu_- + \Delta_{ \mu},  &  \,  t > 0. 
  \end{cases} 
  \label{eq:force3}
\end{equation}
For each path, we combine nonautonomous $(\Delta_\mu,\varepsilon)$ tipping diagrams 
for monotone and nonmonotone shifts to uncover four possible regions:
\begin{itemize}
\item[$\bullet$] {\em Points of tracking} are defined as $(\Delta_\mu,\varepsilon)$ 
settings where the system avoids tipping for monotone and 
nonmonotone shifts. This is the safe region of tracking, sometimes referred to as the ``safe operating space"~\cite{scheffer2015creating}.
\item[$\bullet$] {\em Points of return} are defined as $(\Delta_\mu,\varepsilon)$ 
settings where the system tips for monotone shifts, but does not tip 
for nonmonotone shifts. Here, an otherwise imminent tipping is prevented
by the parameter-shift reversal.
\item[$\bullet$] {\em Points of no return} are defined as $(\Delta_\mu,\varepsilon)$ 
settings where the system tips for monotone and nonmonotone shifts.
Here, tipping is not prevented by the parameter-shift reversal.
\item[$\bullet$] {\em Points of return tipping} are defined as
$(\Delta_\mu,\varepsilon)$ settings where the system does not tip 
for monotone shifts, but tips for nonmonotone shifts. Here, 
the parameter-shift reversal inadvertently induces tipping in an otherwise safe 
situation.
\end{itemize}
Note that the existence, shape, and location of the four regions in the 
$(\Delta_\mu,\varepsilon)$ tipping diagram will, in general, depend on 
the geometric form of the shift $\mu(t)$.

\begin{figure}[t]
\begin{center}
\hspace*{-0.5cm}
\includegraphics[width=31pc]{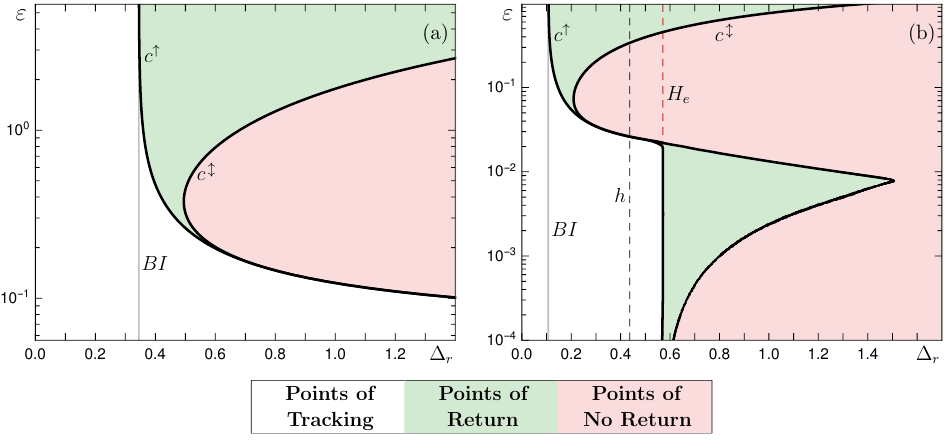}
\end{center}
\caption{
Nonautonomous tipping diagrams in the $(\Delta_r,\varepsilon)$ plane for {\normalfont(a)} R-tipping alone and {\normalfont(b)} R- and B-tipping  are
partitioned into (white) ``points of tracking," (green) ``points of return,"
and (pink) ``points of no return."  
{\normalfont(a)} The nonautonomous tipping-tracking bifurcation curves $c^{\uparrow}$ 
and  $c^{\,\updownarrow}$ are obtained using  monotone  \eqref{eq:force3} and 
 nonmonotone  \eqref{eq:force2} parameter shifts, respectively,  from 
 $p_1 = (1.0,0.075)$ along the parameter path $\Delta_r$ 
 with fixed $m=0.075$ and time-varying $r >1$.
{\normalfont(b)}
The nonautonomous tipping-tracking bifurcation curves $c^{\uparrow}$ 
and  $c^{\,\updownarrow}$ are obtained using  monotone  \eqref{eq:force3} and 
 nonmonotone  \eqref{eq:force2} parameter shifts, respectively,  from 
 $p_1 = (0.5,0.12)$ along the parameter path $\Delta_r$ from Figure~{\normalfont\ref{fig:pht4}(a)}.
 $b=b_c=0.025$.
}
\label{fig:fbs9}
\end{figure} 

\subsection{The ecosystem model}
For the nonautonomous ecosystem model~\eqref{eq:dPdt_na}--\eqref{eq:dHdt_na}, we consider 
two different parameter paths giving rise to two diagrams in 
Figure~\ref{fig:fbs9}. 

The nonautonomous $(\Delta_r,\varepsilon)$ tipping diagram in Figure~\ref{fig:fbs9}(a) 
is obtained for a parameter path with a fixed $m=0.075$, $r_- = 1$, and $r(t)> 1$ 
such that the path crosses the boundary of the basin instability $BI$
but does not cross any autonomous bifurcation. 
Points of no return are bounded by the 
nonautonomous R-tipping bifurcation curve $c^{\,\updownarrow}$ for the 
nonmonotone shift~\eqref{eq:force2}. Points of return are located between 
$c^{\,\updownarrow}$ and the nonautonomous R-tipping bifurcation curve $c^{\uparrow}$ 
 for the monotone shift~\eqref{eq:force3} with $\mu = r$. 
At higher $\varepsilon$, (green) points of return extend over the entire 
$\Delta_r$ interval past the boundary of $BI$. This is indicative of R-tipping 
occurring after the input $r(t)$ reaches its maximum. 
However,
as $\varepsilon$ is decreased, $c^{\uparrow}$ and $c^{\,\updownarrow}$ approach each 
other so that the (green) points of return shrink and appear to vanish at 
$\varepsilon\approx 0.2$. Overlapping of $c^{\uparrow}$ and $c^{\,\updownarrow}$ gives 
rise to apparently direct transitions from (white) tracking to 
(pink) points of no return. This is indicative of R-tipping occurring before the input
$r(t)$ reaches its maximum. 
In other words, the system R-tips from $e_3(t)$ to $e_2(t)$ during the upshift in $r(t)$, 
and the parameter-shift reversal has no effect on the response of the system.

The $(\Delta_r,\varepsilon)$ tipping diagram in Figure~\ref{fig:fbs9}(b) 
is obtained for the parameter path $\Delta_r$ from Figure~\ref{fig:pht4}(a) 
with a fixed $m=0.12$, $r_- = 0.5$, and $r(t)> 0.5$ such that the path 
crosses the boundary  of the basin instability $BI$
as well as the (dangerous) subcritical Hopf bifurcation $H_e$.
At higher $\varepsilon$, and between the $BI$ boundary and $H_e$, where R-tipping is the 
tipping mechanism, the  
diagram is the same as 
in Figure~\ref{fig:fbs9}(a). At intermediate $\varepsilon$, the interplay between 
B-tipping and R-tipping gives rise to a deep wedge in $c^{\,\updownarrow}$, 
which opens up another (green) region with points of return. 
At lower $\varepsilon$, where B-tipping is the tipping mechanism,
the (green) region with points of return shrinks 
but does not seem to vanish as $\varepsilon\to 0$.

\begin{figure}[t]
\begin{center}
\hspace*{-0.5cm}
\includegraphics[width=13cm]{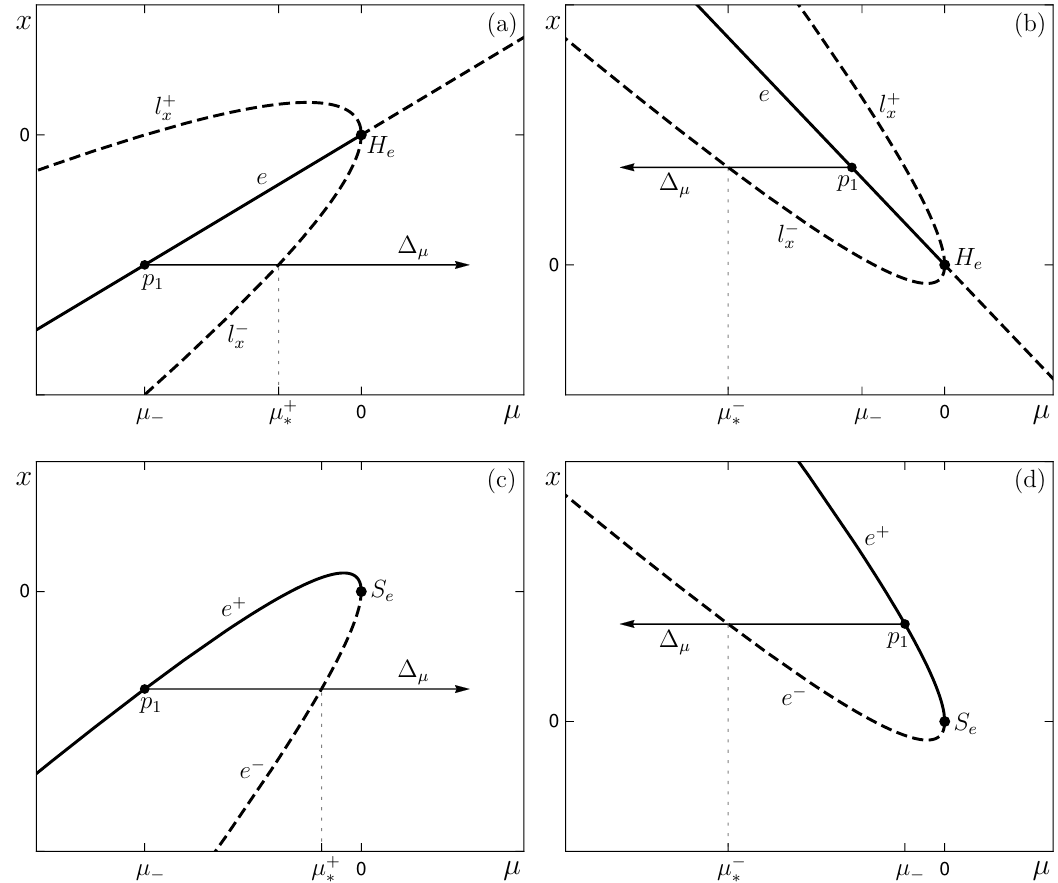}
\end{center}
\caption{
One-parameter bifurcation diagrams for the (tilted) subcritical Hopf normal form 
 \eqref{eq:nf8} with {\normalfont(a)} $s>0$ and {\normalfont(b)} $s<0$, and for the (tilted) saddle-node 
normal form  \eqref{eq:snf1} with {\normalfont(c)} $s>0$ and {\normalfont(d)} $s<0$. Shown are branches of (solid) stable and (dashed) unstable equilibria $e$, branches of the maxima $l_x^+$ and minima $l_x^-$ 
of the $x$-component of the unstable limit cycle,  parameter paths $\Delta_\mu$ from $p_1=\mu_-$, 
and the corresponding boundary $\mu^\pm_*$ of the basin instability {\normalfont(a)}--{\normalfont(b)} $BI(e,\mu_-)$ and {\normalfont(c)}--{\normalfont(d)} 
$BI(e^+,\mu_-)$.
}
\label{fig:S_H_NF1}
\end{figure}

\subsection{The two generic dangerous bifurcations of equilibria}


Here, we obtain nonautonomous tipping diagrams for the saddle-node and subcritical Hopf normal forms to identify typical effects of nonmonotone shifts across a dangerous bifurcation.
While the normal forms are valid close to the bifurcation point, the unstable equilibrium (saddle-node) and unstable limit cycle (Hopf) do capture the global effects of the basin boundary.
Furthermore, we modify the normal forms to mimic the global effects away from the bifurcation point. The main modification involves an additional parameter $s$ that ``tilts'' the branches of solutions in the one-parameter bifurcation diagram giving rise to basin instability; see Figure~\ref{fig:S_H_NF1}. There is also an additional parameter $\alpha$ that quantifies the amount of shear in the Hopf normal form.

\subsubsection{Modified subcritical Hopf normal form}

Consider a system in $\mathbb{R}^2$ akin to the normal form of a subcritical Hopf bifurcation~\cite[sect.3.4]{kuznetsovelements} written in terms of a complex variable $z = x + i y$:
\begin{equation}
\label{eq:nf8} 
\dot z = \left(\mu + i\left[\omega + \alpha \left|z - \mu s \right|^2\right]\right) \left(z - \mu s\right) + \left|z - \mu s \right|^2 (z - \mu s),
\end{equation}
where $\mu$ is the bifurcation parameter, $\omega$ is the angular frequency of small-amplitude oscillations, $\alpha$ quantifies the amount of shear or amplitude-phase coupling, and $s$ is the ``tilt'' parameter. The subcritical Hopf normal form is recovered when we set $s=0$ and  apply a change of 
coordinates to transform away the term proportional to $\alpha$~\cite[sect.3.4]{kuznetsovelements}. There is one 
branch of equilibria,
$$
e(\mu,s) = \mu s + 0i,
$$
that is stable for $\mu<0$ and unstable for $\mu>0$, and one branch of unstable limit cycles,
$$
l(\mu,s,t) = \mu s +\sqrt{-\mu}\,e^{i(\omega - \alpha\mu)t},
$$
that exists for $\mu<0$. The real part of the limit cycle solution oscillates between 
$$
l_x^- (\mu,s)= -\sqrt{-\mu} + \mu s\;\; \mbox{and}\;\; l_x^+(\mu,s) = \sqrt{-\mu} + \mu s,
$$
as shown in Figure~\ref{fig:S_H_NF1}(a)--{\normalfont(b)}.
For every $s\ne 0$, there are two basin instability boundaries:
$$
\mu_*^- = \mu_- -\frac{1 + \sqrt{1-4s^2\mu_- }}{2s^2} < \mu_-
\;\;\mbox{and}\;\; 
\mu_*^+ = \mu_- -\frac{1 - \sqrt{1-4s^2\mu_- }}{2s^2} > \mu_-.
$$

\begin{figure}[t]
\begin{center}
\hspace*{-0.5cm}
\includegraphics[width=31pc]{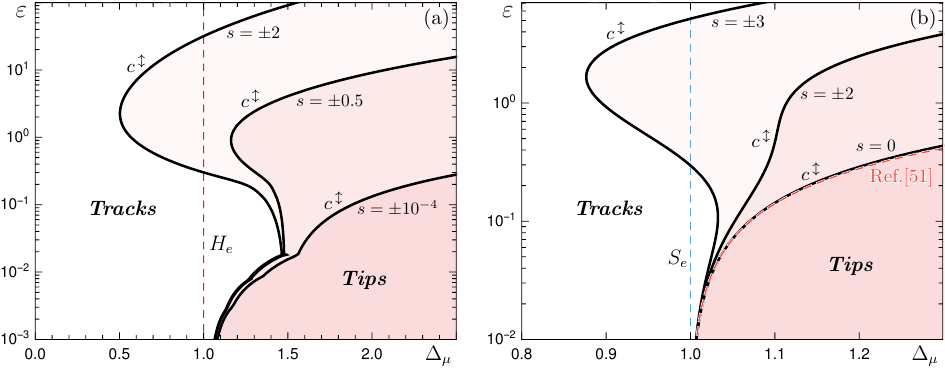}
\end{center}
\caption{
Nonautonomous tipping diagrams in the $(\Delta_\mu,\varepsilon)$ 
 parameter plane for nonmonotone shifts~\eqref{eq:force2} 
 from $p_1=\mu_- =-1$
along the parameter path $\Delta_\mu$ from Figure~\ref{fig:S_H_NF1}{\normalfont(a)} and {\normalfont(c)}. 
{\normalfont(a)} Nonautonomous tipping-tracking bifurcation curves $c^{\,\updownarrow}$
for the (tilted) subcritical Hopf normal form  \eqref{eq:nf10}
with $\alpha=1$, $\omega=1$ and different values of $s$.
{\normalfont(b)}  Nonautonomous tipping-tracking bifurcation curves $c^{\,\updownarrow}$
for the (tilted) saddle-node normal form  \eqref{eq:snf4}
with different values of $s$. The dashed red curve in {\normalfont(b)} is the approximation 
to $c^{\,\updownarrow}$ obtained in  {\normalfont\cite{ritchie2017inverse}} for 
$\varepsilon s$ small enough.
}
\label{fig:HNF2}
\end{figure}

Now, consider the corresponding nonautonomous system
\begin{equation}
\label{eq:nf10} 
\dot z = \left(\mu(t) + i\left[\omega + \alpha \left|z - s\,\mu(t)  \right|^2\right]\right) \left(z - s\,\mu(t)\right) + \left|z - s\,\mu(t) \right|^2 (z - s\,\mu(t)).
\end{equation}
First, we analyze R-tipping for nonmonotone $\mu(t)$ given by  \eqref{eq:force2}
where we replace $r$ with $\mu$ and use $\mu_- = -1$ and  $\Delta_\mu >0$ (Figure~\ref{fig:HNF2}(a)).
Tipping from the stable equilibrium $e$ requires nonzero 
$s$ because the branch of equilibria $e=\mu s + 0i$ is flow-invariant
when $s=0$.
For $s=10^{-4}$, we obtain  $\mu_*^+\approx -10^{-8}$, meaning that the region of 
basin instability between $\mu_*^+$ and $H_e$ is negligible. The only tipping 
that occurs in the nonautonomous system 
is  B-tipping for $\Delta_\mu > 1$, as evidenced by the tipping-tracking 
transition curve $c^{\,\updownarrow}$ in the $(\Delta_\mu,\varepsilon)$ tipping diagram.
When $s=0.5$, the basin instability boundary moves to $\mu_*^+=2\sqrt{2}-3\approx-0.17$
or $\Delta_\mu \approx 0.83$, and the region of basin instability becomes nonnegligible. 
As a result, the  curve $c^{\,\updownarrow}$ deviates from the case $s=10^{-4}$ in different ways.
While R-tipping still does not occur, basin instability gives rise to
a tongue/fold on $c^{\,\updownarrow}$ and a range of shift magnitudes $\Delta_\mu$ 
with three critical rates. 
When the ``tilt'' is increased to $s=2$, the basin instability boundary moves to 
$\mu_*^+\approx-0.61$ or $\Delta_\mu \approx 0.39$. Now, in addition to B-tipping and a range of $\Delta_\mu$ 
with three critical rates, there is R-tipping for $\Delta_\mu<1$. The tracking-tipping transition curve $c^{\,\updownarrow}$ with the R-tipping tongue at higher rates 
and the ``slow'' approach toward $H_e$ as $\varepsilon\to 0$
closely resembles the tipping diagram for the ecosystem model from Figure~\ref{fig:pht4}(b). 
The most noticeable difference from the ecosystem model 
is the absence of the deep wedge at the intermediate rates, possibly due to 
the absence of the homoclinic bifurcation $h$. Instead, there is 
a characteristic kink on the $c^{\,\updownarrow}$ curves  near $\varepsilon=10^{-2}$ 
in Figure~\ref{fig:HNF2}(a)  with multiple wiggles such as those shown in the inset
of Figure~\ref{fig:fbs10}(b). The origin of the kink and the wiggles, as well as the scaling 
law for  $c^{\,\updownarrow}$  in the limit $\varepsilon\to 0$, may be related to 
so-called buffer points~\cite{neishtadt1987persistence,neishtadt1988persistence} 
and are left for future study.

\overfullrule=0pt
The same is true for ``points of return" and ``points of no return" shown in Figure~\ref{fig:SNHNF}(b1), 
where the tracking-tipping transition curve $c^{\uparrow}$ is obtained for the monotone 
parameter shift~\eqref{eq:force3}. Interestingly, for a sufficiently high ``tilt'' 
parameter $s$, a new region of ``points of return tipping" appears in the 
diagram (Figure~\ref{fig:SNHNF}(c1)) that is not present in the ecosystem model. 
This means that, in general, all four regions identified in the beginning of 
section~\ref{sec:pnr} can be present for a nonmonotone passage through a subcritical 
Hopf bifurcation. Furthermore, the rotational symmetry in the phase space of the (modified) 
Hopf normal form implies a symmetry in the basin instability boundaries
\begin{equation}
\label{eq:hopftilt}
\mu_*^\pm(s) = \mu_*^\pm(-s),
\end{equation}
meaning that the system has the same basin instability properties for 
$s$ and $-s$. Thus, we obtain the same tipping diagrams for $s$ and $-s$ in the left column of Figure~\ref{fig:SNHNF}, in line with our  R-tipping criterion from section~\ref{sec:tcTtip}. 
Finally, for a fixed $s\ne 0$, R-tipping for an increasing $\mu(t)$ requires a smaller shift 
magnitude than R-tipping for the decreasing $\mu(-t)$. This is why the region of 
``points of return tipping'' in Figure~\ref{fig:SNHNF}(c1) is small.

\subsubsection{Modified saddle-node normal form}

Consider a system in $\mathbb{R}$ akin to the normal form of a saddle-node 
bifurcation~\cite[sect.3.2]{kuznetsovelements}:
\begin{equation}
\label{eq:snf1} 
\dot x = -(x-\mu s)^2 -\mu,
\end{equation}
where $\mu$ is the bifurcation parameter and $s$ is the ``tilt'' parameter. The branches of stable $e^+$ and unstable $e^-$ equilibria exist for $\mu\le 0$ and are given by
\begin{equation}
\label{eq:snfeqsol1} 
e^+ (\mu,s)= \mu \, s + \sqrt{-\mu}, \;\;\mbox{and}\;\; e^-(\mu,s)= \mu \, s -\sqrt{-\mu},
\end{equation}  
as shown in Figure~\ref{fig:S_H_NF1}(c)--(d). The basin instability boundary is given by
\begin{equation}
\label{eq:snf3} 
\mu_* = - \left(\sqrt{-\mu_-} - \frac{1}{s}\right)^2.
\end{equation}

\begin{figure}[t]
\begin{center}
\hspace*{-0.5cm}
\includegraphics[width=31pc]{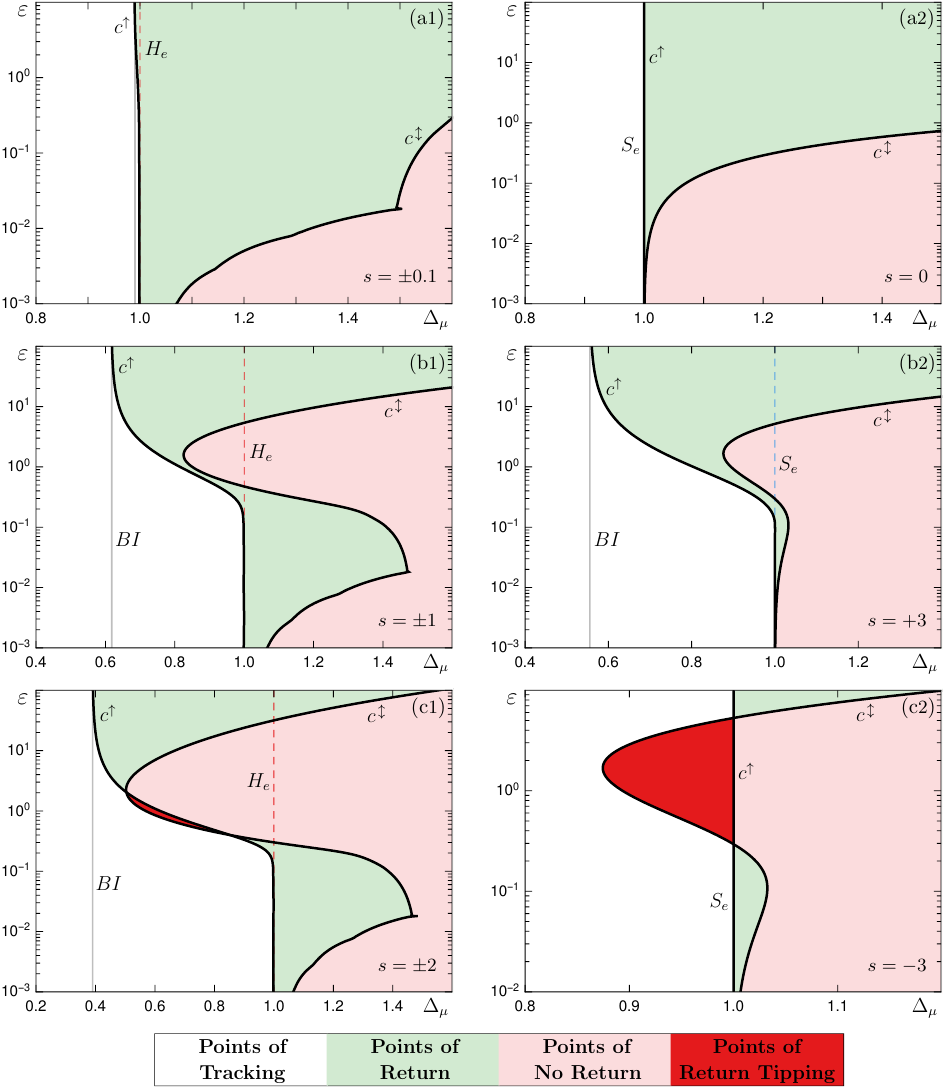}
\end{center}
\caption{
Nonautonomous tipping diagrams  for {\normalfont(a1)--(c1)} the (tilted) Hopf normal form~\eqref{eq:nf10}
 and {\normalfont(a2)--(c2)} the (tilted) saddle-node normal form~\eqref{eq:nf10} with different values of 
 $s$ are partitioned into (white) ``points of tracking," (green) ``points of return," (pink) 
 ``points of no return" and (red) ``points of return tipping." 
The nonautonomous tipping-tracking bifurcation curves $c^{\uparrow}$ and  $c^{\,\updownarrow}$ are obtained
for monotone  \eqref{eq:force3} and nonmonotone  \eqref{eq:force2} 
parameter shifts with  $\mu_-=-1$ along the parameter paths $\Delta_\mu$ from Figure~\ref{fig:S_H_NF1}{\normalfont(a)} and {\normalfont(c)}, respectively.
}
\label{fig:SNHNF}
\end{figure} 

Now, consider the corresponding nonautonomous system 
\begin{align}
\label{eq:snf4} 
\dot x = -(x-\mu(t) s)^2 -\mu(t).
\end{align}
First, we analyze R-tipping for nonmonotone $\mu(t)$ given by  \eqref{eq:force2}
where we replace $r$ with $\mu$ and use $\mu_- = -1$ and  $\Delta_\mu >0$ (Figure~\ref{fig:HNF2}(b)). 
When $s=0$, there is no basin instability and R-tipping cannot occur~\cite[Theorem 3.2(1)]{ashwin2017parameter}.
The only tipping that occurs for $s=0$ is B-tipping for $\Delta_\mu>1$. 
The tracking-tipping transition curve $c^{\,\updownarrow}$ in the $(\Delta_\mu,\varepsilon)$ 
tipping diagram is in very good agreement with the critical ``exceedance time" inverse square law
$$
t_e \approx \frac{2}{\sqrt{\Delta_\mu + \mu_-}},
$$
derived 
in  \cite{ritchie2017inverse} for $\varepsilon s$ small enough.
To demonstrate the agreement, we use  \eqref{eq:force2} to rewrite 
the $t_e$ formula above in terms of $\varepsilon$ and $\Delta_\mu$,
\begin{equation}
\label{eq:texceedtrans}
\varepsilon \approx \sqrt{\Delta_\mu +\mu_-}\;\mbox{sech}^{-1}\left(\frac{-\mu_-}{\Delta_\mu} \right),
\end{equation}
and plot condition~\eqref{eq:texceedtrans} as a dashed red curve in Figure~\ref{fig:HNF2}(b); 
see~\cite[sect.6.2.2]{OKeeffePhD2020} for more details. 
However, for nonzero $s$ the tracking-tipping transition curve $c^{\,\updownarrow}$
can deviate from the inverse square law noticeably and qualitatively. When $s=2$, 
the tracking-tipping transition curve $c^{\,\updownarrow}$ deviates 
from the case $s=0$ noticeably (up to an order of magnitude in $\varepsilon$), 
but the changes are quantitative, and R-tipping does not occur 
despite crossing the basin instability boundary at $\mu_*=-1/4$. 
When the ``tilt'' is increased to $s=3$, the basin instability boundary moves to 
$\mu_*=-4/9$, meaning that $e^+$ is basin unstable for $\Delta_\mu > \mu_*-\mu_-=5/9$. 
This results in two significant changes to the tracking-tipping transition curve $c^{\,\updownarrow}$. 
First, $c^{\,\updownarrow}$  develops two folds and becomes S-shaped, giving rise to 
a range of shift magnitudes $\Delta_\mu$ with three different critical rates. 
Second, in addition to B-tipping, there is an R-tipping tongue for $\Delta_\mu<1$.
In contrast to the subcritical Hopf bifurcation, the $c^{\,\updownarrow}$ curves 
clearly converge to $S_e$ as $\varepsilon\to 0$. This is because,
as $\varepsilon\to 0$, the solution 
jumps off the branch of stable equilibria at the bifurcation point with
no delay, meaning there is no time to turn around and avoid tipping~\cite{berglund2006noise,majumdar2013transitions}.

Analysis of  ``points of return" and ``points of no return" near a saddle-node 
bifurcation reveals much similarity to the subcritical 
Hopf bifurcation when $s>0$, but not when $s<0$ (Figure~\ref{fig:SNHNF}(b2)). 
The striking difference for 
$s=-3$ is the large region of ``points of return tipping,'' where  there is 
R-tipping for nonmonotone $\mu(t)$, but not for monotone increasing $\mu(t)$
(Figure~\ref{fig:SNHNF}(c2)). This difference is a consequence 
of asymmetry in the (modified) saddle-node normal form. To be more specific, 
\begin{equation}
\label{eq:sntilt}
\mu_*(s) \ne \mu_*(-s),  
\end{equation}
meaning that the system has different basin instability properties for 
$s$ and $-s$.  According to the  R-tipping criterion from section~\ref{sec:tcTtip},
given a suitable $\mu(t)$ that increases over time,
the system is guaranteed to R-tip for $s>0$, but not for $s<0$. Conversely, 
given a suitable $\mu(t)$ that decreases over time, the system is guaranteed
to R-tip for $s<0$, but not for $s>0$. Thus, ``points of return tipping" 
cannot occur for $s>0$ and are expected to occur for $s<0$, which 
explains the diagrams for $s=3$ and $s=-3$ in Figure~\ref{fig:SNHNF}(b2) and (c2).

\subsubsection{Universal properties of tipping near a dangerous bifurcation}

A comparison between the tracking-tipping transition curves $c^{\,\updownarrow}$ for 
the modified subcritical Hopf (Figure~\ref{fig:HNF2}(a)) and saddle-node (Figure~\ref{fig:HNF2}(b)) 
normal forms reveals some universal tipping properties. In both systems, the tracking-tipping 
transition curve $c^{\,\updownarrow}$ becomes S-shaped, gives rise to three critical rates, 
and develops an R-tipping tongue as the ``tilt'' parameter $s$ is increased. 

On the other hand, multiple critical rates and R-tipping are achieved for a smaller 
`tilt' parameter $s$ in the modified Hopf normal form, whereas the approach of $c^{\,\updownarrow}$ 
toward the bifurcation as $\varepsilon\to 0$ is much faster and follows a different scaling 
law in the modified saddle-node normal form. Finally, owing to the
different basin instability properties for $s$ and $-s$, 
as shown by Eqs.~\eqref{eq:hopftilt} and~\eqref{eq:sntilt}, 
a saddle-node bifurcation may give rise to a larger region of ``points of return tipping."

\section{Conclusion}
\label{sec:concl}

We analyze nonlinear tipping phenomena in systems with time-varying
inputs, using examples of an ecological model~\cite{scheffer2008pulse} and modified 
saddle-node and subcritical Hopf normal forms with parameter shifts. 
The ecological model exhibits a somewhat counterintuitive behavior, where the herbivore
population persists for a slow increase in the food growth rate but tips to extinction 
when the food growth rate increases too fast.
We analyze such tipping phenomena as nonautonomous bifurcations. The
proposed mathematical framework uses the global property  of  basin  
boundaries  in  the  autonomous  frozen  system  
with fixed-in-time inputs to give criteria for the 
occurrence of nonautonomous bifurcations in the system with 
time-varying inputs. This framework aims to be easily accessible 
to applied scientists, addressing two questions of relevance:  
critical factors for tipping and the possibility of preventing tipping 
by a trend reversal. Our results give new insight into the sensitivity of 
ecosystems to the magnitudes and rates of environmental change. 

Genuine nonautonomous {\em rate-induced bifurcations (R-tipping)}, 
which are entirely due to the rate of change of the input parameters, 
are shown to correspond to {\em maximal canard trajectories} that 
track moving unstable states for an infinite time. We give simple 
criteria for the occurrence of R-tipping in the nonautonomous 
system using the concepts of {\em parameter paths} and 
{\em basin instability on a path} in the autonomous frozen 
system. These criteria will allow applied scientists to easily 
test whether their systems have critical rates of change 
and uncover new phenomena, such as {\em R-tipping tongues},
in the  nonautonomous tipping diagrams. 
We also note that R-tipping problems can be transformed into connecting (heteroclinic) orbit problems using a suitable compactification technique developed in  \cite{wieczorek2020}.

Classical autonomous bifurcation analysis of the frozen system reveals 
a codimension-three degenerate Bogdanov--Takens bifurcation. This is the source 
of a (dangerous) subcritical Hopf bifurcation and is an organizing center for {\em bifurcation-induced tipping (B-tipping)} in the nonautonomous ecosystem model.
Superimposing regions of basin instability onto classical bifurcation diagrams 
adds information about genuine nonautonomous bifurcations, which
can be even more relevant in certain applications, but are missed by classical autonomous stability analysis.
Thus, our approach gives a comprehensive insight into system stability, beyond classical autonomous bifurcations and adiabatic effects of a parameter change. 


Analysis of the interaction between B-tipping and  R-tipping
reveals an  {\em S-shaped nonautonomous bifurcation curve}  with multiple 
critical rates in the nonautonomous tipping diagram. 
This curve captures different tipping mechanisms,
giving rise to
{\em points of tracking}, {\em points of return} where tipping can be prevented
by the parameter-trend reversal, {\em points of no return} where tipping cannot be 
prevented by the reversal, and {\em points of return tipping} where tipping is 
inadvertently induced by the reversal. Analysis of the modified  
saddle-node and subcritical Hopf normal forms suggests that these features could be considered universal  for nonmonotone parameter shifts that cross a basin instability boundary and a dangerous bifurcation, and then reverses.

\section*{Acknowledgment}
The first author  thanks M. Mortell for useful discussions on singular perturbations.

\bibliographystyle{siamplain}
\bibliography{references}
\end{document}

%% file: ex_shared.tex
\usepackage{amsfonts}
\usepackage{graphicx}
\usepackage{epstopdf}
\usepackage{algorithmic}
\usepackage{graphics}
\usepackage{amsmath}
\usepackage{amssymb}
\usepackage{verbatim}
\usepackage[utf8]{inputenc}
\usepackage[english]{babel}

\usepackage{tikz} 
\usepackage{listings}
\usepackage{wrapfig}
\usepackage[abs]{overpic}
\usepackage{pict2e}
\usepackage{float}
\usepackage{cite}
\usepackage{hyperref}
\newcommand*\circled[1]{\tikz[baseline=(char.base)]{
\node[shape=circle,draw,inner sep=1](char) {#1};}}

\definecolor{ggray}{HTML}{D8D8D8}

\newcommand{\eps}{\varepsilon}

\newsiamremark{remark}{Remark}
\newsiamremark{hypothesis}{Hypothesis}
\crefname{hypothesis}{Hypothesis}{Hypotheses}
\newsiamthm{claim}{Claim}

\headers{Tipping Phenomena and Points of No Return}{P.\,E. O'Keeffe, S. Wieczorek}

\title{Tipping Phenomena and Points of No Return in Ecosystems: Beyond Classical Bifurcations}
\author{{Paul E. O'Keeffe and Sebastian Wieczorek}\thanks{Department of Applied Mathematics, University College Cork, Ireland 
  (\email{paulokeeffe@umail.ucc.ie}, \email{sebastian.wieczorek@ucc.ie}).}
  \funding{The work of the second author was supported by the CRITICS
Innovative Training Network via the European Union’s Horizon 2020 research and innovation programme under grant agreement 643073.}}

\usepackage{amsopn}
